\documentclass[11pt]{article}

\usepackage[margin=1in]{geometry}
\usepackage{amsmath,amssymb,amsthm,bm}
\usepackage{booktabs}
\usepackage{graphicx}
\usepackage{float}
\usepackage{microtype}
\usepackage[authoryear,round]{natbib}
\usepackage{algorithm}
\usepackage{algpseudocode}
\usepackage[hidelinks]{hyperref}
\usepackage{placeins}
\usepackage{cleveref}

\newcommand{\E}{\mathbb{E}}
\newcommand{\R}{\mathbb{R}}
\newcommand{\indep}{\mathrel{\perp\!\!\!\perp}}
\newcommand{\cF}{\mathcal{F}}
\newcommand{\cG}{\mathcal{G}}

\title{The Monge Optimal Transport Barycenter Problem\\
with Adaptive Features}
\author{Yige Zhou\thanks{New York University Shanghai.}
\and
Esteban G. Tabak\thanks{Courant Institute of Mathematical Sciences,
New York University.}}
\date{}

\begin{document}
\maketitle

\begin{abstract}
The data-driven Monge optimal transport barycenter problem seeks a map $y=T(x,z)$ that minimally perturbs an outcome $x$ while removing its statistical dependence on covariates $z$. Numerical implementations approximate the independence constraint through finite-dimensional test-function spaces, making the effectiveness of dependence removal sensitive to the choice of features. 
This article develops a procedure whereby the features themselves are optimized throughout the procedure. In particular, we adopt a feature space for $f$ consisting of Gaussian radial basis functions (RBFs) with adaptive individual centers and bandwidths.
Synthetic experiments show that while fixed nonlinear RBF features can capture dependence structures that are poorly represented by linear features, adapting the RBF centers and bandwidths can further improve dependence removal within a limited feature space and achieve residual dependence comparable to that obtained with a larger fixed basis. We further extend the formulation to multidimensional covariates and evaluate its applicability to real-world problems, including removing age-related dependence in biomedical measurements and dependence on temperature and humidity in bike-sharing demand.
\end{abstract}

\noindent\textit{Keywords:}
Optimal transport; Statistical independence;
Adaptive feature representations; Gaussian radial basis functions;
Dependence removal; Numerical optimization

\section{Introduction}

Let $(x,z)\sim\pi$ denote the joint distribution of outcome variables $x\in\mathcal X\subseteq\R^{d_x}$ and observed covariates $z\in\mathcal Z\subseteq\R^{d_z}$.
The Monge optimal transport barycenter problem (OTBP) seeks a
$z$-dependent map that removes all dependence of $x$ on $z$ with minimum transportation cost:
\begin{equation}
  \min_{T}\;\E_{\pi}\!\left[c\bigl(x,T(x,z)\bigr)\right]
  \quad\text{subject to}\quad
  y=T(x,z)\indep z .
  \label{eq:otbp}
\end{equation}
For the canonical cost $c(x,y)=\tfrac12\lVert x-y\rVert^2$, the distribution
of $y$ is the Wasserstein barycenter of the family
$\{\rho(\cdot\mid z)\}$ whenever the optimal couplings are supported on maps
\citep{agueh2011barycenter,pass2013infinitely,tabak2018variability,tabak2022distributional,lipnick2026monge}.

The independence constraint in \cref{eq:otbp} is equivalent to the uncorrelation between all measurable functions $f(z)$ and $g(y)$:
\begin{equation}
  \E[f(z)g(y)]-\E[f(z)]\E[g(y)]=0
  \quad\text{for all measurable }f\text{ and }g.
  \label{eq:independence}
\end{equation}
When only finite samples of $\pi$ are available, one must restrict the allowable functions $f$ and $g$ to those lying in suitable finite-dimensional spaces $\cF$ and $\cG$ \citep{gretton2005measuring}. This creates a statistical approximation problem, since any dependence invisible to $\cF\times\cG$ remains in the transported variable $y$.

In this article, we focus on feature representations with optimizable parameters, in contrast to conventional fixed feature representations. We begin with the Monge formulation and the test-function characterization of independence, followed by a two-step numerical procedure. The first step identifies statistical dependence between the transported variable and the covariates through test functions, while the second updates the transported samples to reduce the detected dependence while controlling transportation cost.

We first consider the simplest setting, with fixed linear features and then, to capture more complex dependence structures, Gaussian radial basis functions $f(z,\alpha)$, where $\alpha$ contains the centers and bandwidths. While a fixed RBF basis provides greater flexibility than linear features, its effectiveness still depends on the choice of these parameters. We therefore allow the centers and bandwidths to be optimized through the dependence objective in the first step and derive the corresponding update formulas, with the aim of improving dependence detection within a limited feature space.

\section{Adversarial and Finite-Dimensional OTBP Formulation} 
In this section, we briefly describe a modified version of the adversarial formulation and finite-dimensional numerical procedure for the Monge optimal transport barycenter problem in \cite{lipnick2026monge}.
\subsection{Adversarial characterization of independence}
The independence condition in (2) can be enforced through test functions.
After centering $f$, independence requires that
\[
\mathbb{E}[f(z)g(y)] = 0
\]
for every admissible pair $f$ and $g$. 
This motivates an
adversarial formulation in which the test functions reveal any
remaining dependence, while the transport map modifies $y$ to remove it.

After separating the amplitudes of the test functions from their shapes and
absorbing the amplitudes into a penalty parameter $\lambda>0$, the population
problem is written as
\begin{equation}
\min_{y=T(x,z)}
\max_{\substack{f\in\mathcal F,\; g\in\mathcal G\\
\mathbb E[f]=0,\;\|f\|=\|g\|=1}}
\left\{
\mathbb E[c(x,y)]
+
\lambda\,\mathbb E[f(z)g(y)]
\right\},
\label{eq:adversarial_otbp}
\end{equation}
%
%
%
where the parameter $\lambda$ controls the trade-off between
removing dependence and keeping the transport cost small.

\subsection{Finite-dimensional numerical implementation}

The unrestricted classes of test functions introduced above are infinite-dimensional, which would overfit any distribution $\pi$ only known through a finite sample.
We therefore restrict these functions to finite-dimensional feature spaces. 
Given observations
$\{(x_i,z_i)\}_{i=1}^n$ and current transported samples
$\{y_i\}_{i=1}^n$, let $F$ and $G$ denote the centered
feature matrices:
\begin{equation}
F_{ij}
=
F^j(z_i)
-
\frac{1}{n}\sum_{r=1}^n F^j(z_r),
\qquad
G_{ik}
=
G^k(y_i)
-
\frac{1}{n}\sum_{r=1}^n G^k(y_r),
\label{eq:centered_features}
\end{equation}
where the $F^j(z)$ and $G^k(y)$ are given functions.
The finite-dimensional test functions are then represented as
\begin{equation}
f(z) = F(z)a,
\qquad
g(y) = G(y)b,
\label{eq:finite_test_functions}
\end{equation}
where $a$ and $b$ are coefficient vectors.
For a finite sample, because $f$ and $g$ are centered, their empirical covariance is given by
\[
\mathbb{E}[f(z)g(y)]
\approx
\frac{1}{n}\sum_{i=1}^n f(z_i)g(y_i)
=
\frac{1}{n}f^\top g.
\]

In order to turn this covariance into a correlation, the proposal in \cite{lipnick2026monge} orthogonalizes the matrices $F$ and $G$ and adds the constraint that $\|a\| = \|b\| =1$. Instead, in the numerical implementation considered below, the explicit unit-norm constraints are replaced by the maximization over $f$ and $g$ of the functional
\begin{equation} 
J(f,g) = \frac{1}{n}f^\top g - \frac{1}{2n^2}\lVert f\rVert^2\lVert g\rVert^2, \label{eq:J} \end{equation} 
%
which selects the same dependence directions as maximizing
the empirical correlation.

The finite-dimensional numerical procedure alternates between two
operations. The first searches for dependence between the current transported
samples and the covariates. The second updates the transported samples, initialized at $y^{(0)} = x$, to
reduce the detected dependence while controlling the transportation cost.

\paragraph{Step 1: find dependence.}
For the current transported samples $y^{(k)}$, the coefficient vectors are
chosen through
\begin{equation}
(a^{(k)},b^{(k)})
\in
\arg\max_{a,b}
J\!\left(
Fa,
G(y^{(k)})b
\right).
\label{eq:step1}
\end{equation}
This step selects the strongest remaining dependence direction detectable by
the chosen feature spaces. 

\paragraph{Step 2: remove dependence.}
Once one or more dependence directions have been identified, they are held
fixed while the transported samples are updated. For a fixed dependence direction $f$ identified in Step~1, define the
corresponding Step~2 penalty by
\[
P(y)
=
\max_b
J\!\left(
f,G(y)b
\right),
\]
where $f$ is held fixed and $b$ is re-optimized for the current transported
samples $y$.
Let $P^{(\ell)}(y)$ denote the dependence penalty associated with
the $\ell$th direction. The Step~2 objective has the form
\begin{equation}
L(y)
=
\frac{1}{2n}
\sum_{i=1}^n
\|y_i-x_i\|^2
+
\sum_{\ell=1}^q
\lambda_\ell P^{(\ell)}(y).
\label{eq:step2_objective}
\end{equation}
The first term is the quadratic transportation cost, while the second term
penalizes the dependence detected in Step~1. Here, $\lambda_\ell$ is an
adaptive, direction-specific penalty weight rather than a fixed tuning
parameter.\footnote{In the MATLAB implementation, the weights are recomputed
at every Step~2 iteration. The current correlations are rescaled as
$\widetilde{\sigma}_\ell=\min\{\max(\sigma_\ell/\sigma_0-0.9,0),20\}$,
where $\sigma_0=\max\{0.001,n^{-1/2}\}$. After sorting the directions, the
algorithm sets $\lambda_\ell=\max\{\widetilde{\sigma}_\ell-c_\ell,0\}$,
where $c_\ell\leq0$ accounts for the alignment between the current accumulated
Step~2 gradient and $\nabla_yP^{(\ell)}(y)$. Thus, the weights are updated
dynamically rather than selected as fixed hyperparameters.}

The gradient used in the implementation is
\[
\nabla_y L
=
\frac{y-x}{n}
+
\sum_{\ell=1}^{q}
\lambda_\ell \nabla_y P^{(\ell)}(y),
\qquad
y \leftarrow y-\eta \nabla_y L.
\]
Thus, Step~1 searches for dependence, whereas Step~2 modifies the transported
samples to reduce the detected dependence while keeping the quadratic
transportation cost small. 

The global algorithm alternates between steps 1 and 2, with step 1 selecting one new test function $f$ each time, and step 2 using all the $\{f_l\}$ selected so far.

\section{Monomials and Nonlinear Features with a fixed RBF basis }
\label{sec:linear}

The simplest linear setting uses fixed linear functions of $z$ and low-order polynomials of $y$. For scalar variables, a representative choice is
\begin{equation}
  F(z)=[z],\qquad G(y)=[y,y^2],
  \label{eq:linear-features}
\end{equation}
followed by column centering.  

This linear setting is simple, inexpensive, and easy to interpret. However,
dependence in real data is often more complex than what can be represented by
linear or low-order polynomial features. Because the dependence condition is approximated through finite-dimensional
feature spaces, the choice of basis functions can be changed without altering
the underlying OTBP formulation. We therefore replace the linear features by
more flexible nonlinear ones to capture a broader class of dependence structures.
The next section develops this nonlinear setting.


\subsection{Nonlinear feature matrix for f}

To extend the setting considered above, we replace the fixed linear feature matrix by a nonlinear feature matrix $F$, constructed from
Gaussian radial basis functions (RBFs) of $z$ \citep{micchelli1986interpolation}.

For one-dimensional $z$, the raw nonlinear feature matrix is
\begin{equation}
F^{\rm raw}_{ij}=\exp\!\left[-\frac{(z_i-z_j^c)^2}{2h^2}\right],
\qquad i=1,\ldots,n,\quad j=1,\ldots,n_c,
\label{eq:fixed-rbf}
\end{equation}
where $z_j^c$ is the $j$th center and $h$ is the bandwidth.  The feature matrix is
obtained by centering every column,
\begin{equation}
F_{ij}=F^{\rm raw}_{ij}-\frac{1}{n}\sum_{r=1}^n F^{\rm raw}_{rj}.
\qquad f=Fa.
\label{eq:centered-rbf}
\end{equation}
%
We initialize the
centers by $k$-means \citep{lloyd1982least} and use Silverman's rule of thumb \citep{silverman1986density} for the shared
bandwidth,
\begin{equation}
h=0.9\min\!\left\{\widehat\sigma_z,
\frac{\operatorname{IQR}(z)}{1.34}\right\}n^{-1/5}.
\label{eq:silverman}
\end{equation}

In the version with fixed nonlinear features, the centers and bandwidth are chosen at the outset and the centered RBF matrix $F$ is then held fixed. 

\subsection{Examples with fixed linear and nonlinear features}

In the synthetic examples below, the data are generated as a
$z$-dependent signal plus an independent noise component. Therefore, in addition to plotting x and y versus z, we also plot the true noise against the transported
variable $y$. This provides a visual diagnostic of whether the dependence
has been removed while keeping the transport-induced distortion small.

\paragraph{Nonlinear heteroscedastic example.}
We first consider
\begin{equation}
x=z+\sigma(z)\epsilon,\qquad
\sigma(z)=0.1\bigl(1+0.8\sin(2\pi z)\bigr),
\qquad \epsilon\sim\mathcal{N}(0,1).
\label{eq:nonlinear-heteroscedastic}
\end{equation}
Here the conditional mean is linear in $z$, but the conditional variance
changes nonlinearly.

\begin{figure}[H]
\centering
\begin{minipage}{0.48\linewidth}\centering\includegraphics[width=0.55\linewidth]{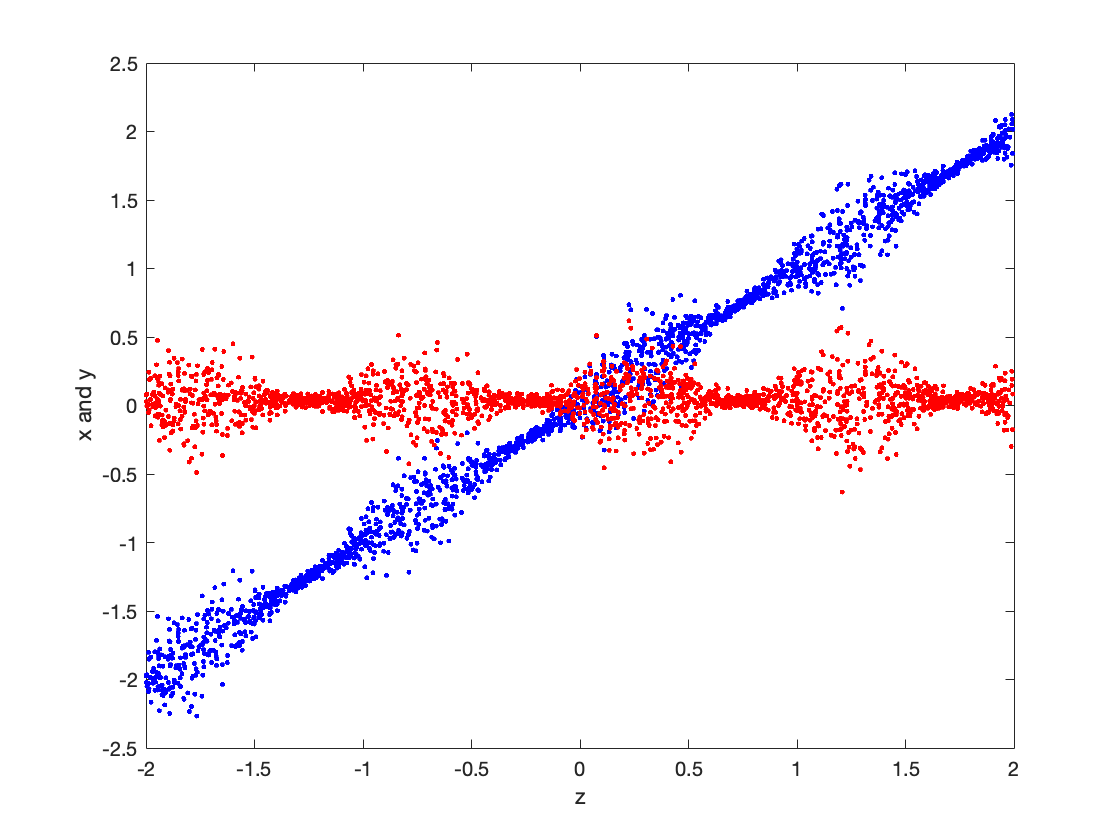}\par
\footnotesize (a) Linear $F$: $x$ and $y$ versus $z$
\end{minipage}\hfill
\begin{minipage}{0.48\linewidth}\centering\includegraphics[width=0.55\linewidth]{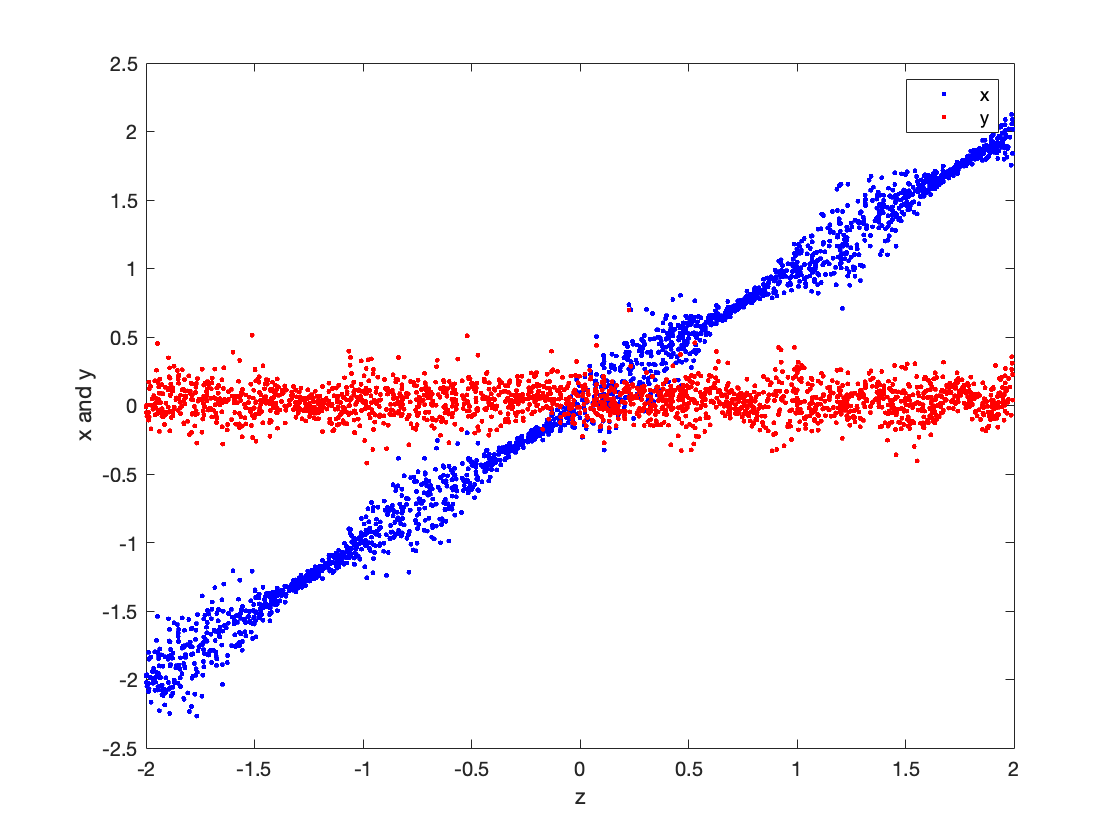}\par
\footnotesize (b) Fixed nonlinear $F$: $x$ and $y$ versus $z$
\end{minipage}

\vspace{0.4em}
\begin{minipage}{0.48\linewidth}\centering\includegraphics[width=0.55\linewidth]{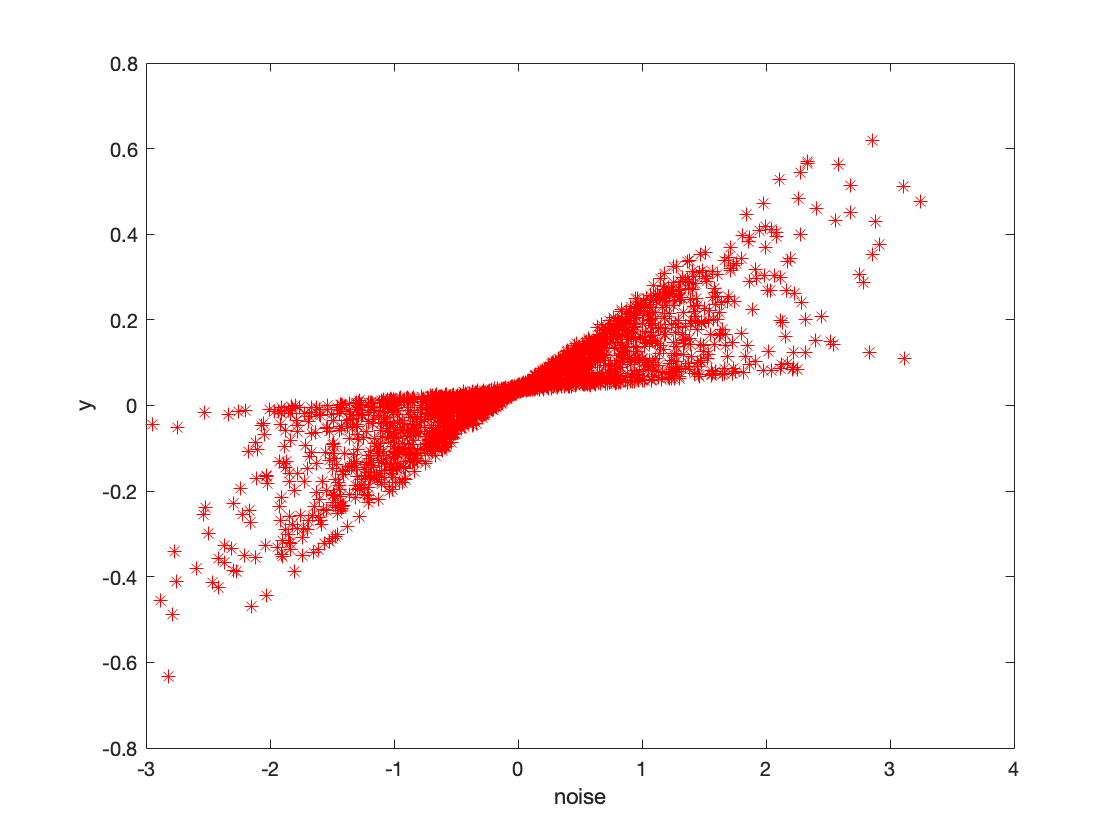}\par
\footnotesize (c) Linear $F$: noise versus $y$
\end{minipage}\hfill
\begin{minipage}{0.48\linewidth}\centering\includegraphics[width=0.55\linewidth]{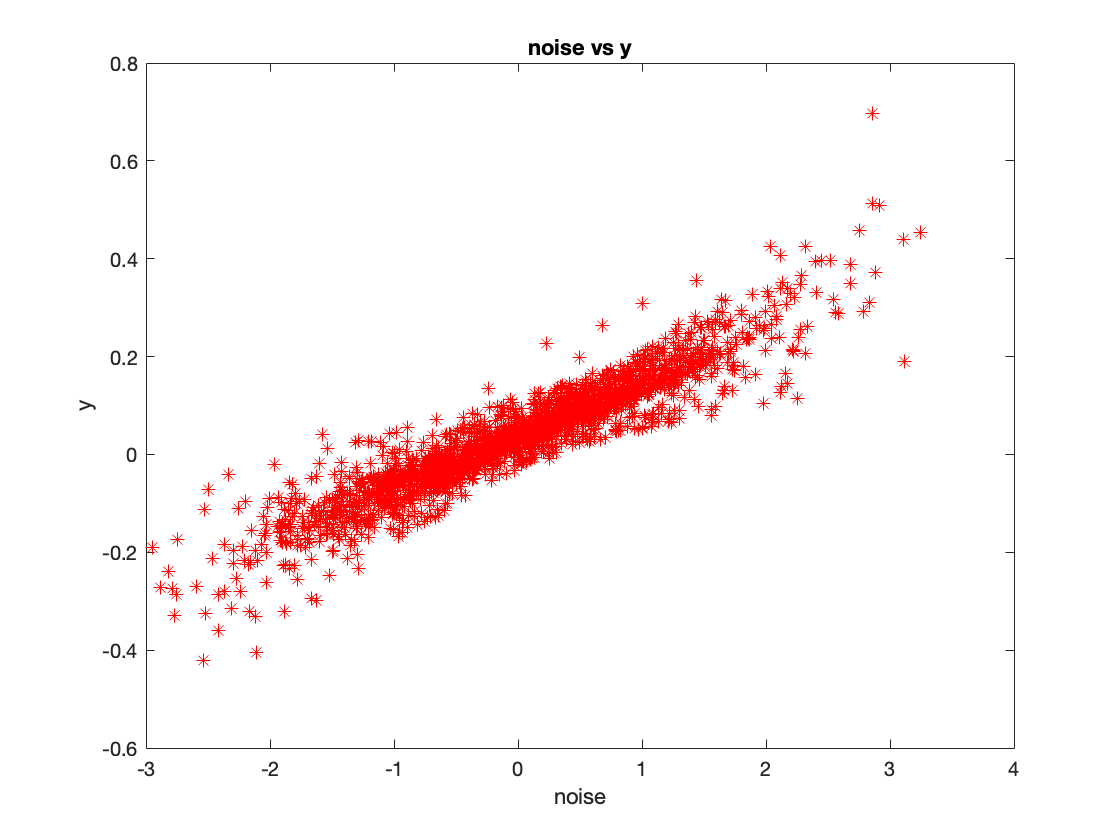}\par
\footnotesize (d) Fixed nonlinear $F$: noise versus $y$
\end{minipage}
\caption{Nonlinear heteroscedastic example. The linear feature matrix removes
the mean dependence but leaves a visible variance pattern. The fixed nonlinear
feature matrix removes both the mean and nonlinear variance dependence.}
\label{fig:nonlinear-heteroscedastic}
\end{figure}

\begin{figure}[H]
\centering
\begin{minipage}{0.48\linewidth}
\centering\includegraphics[width=0.55\linewidth]{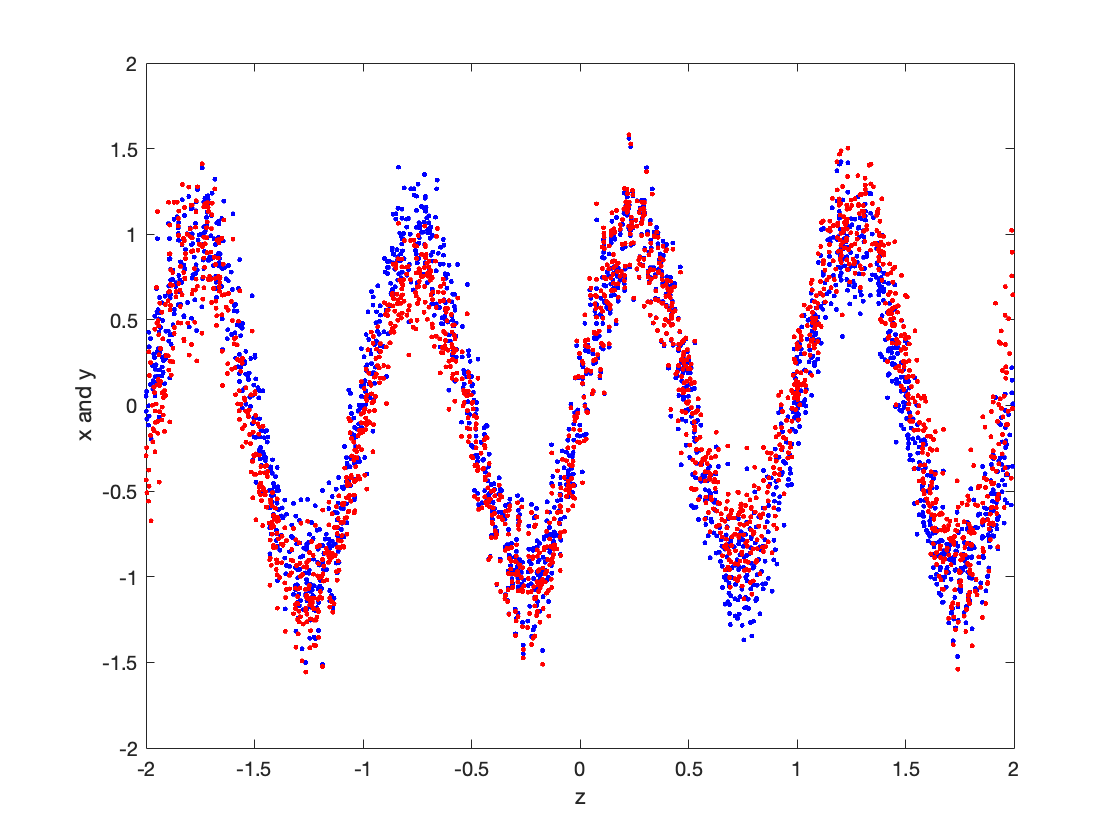}\par
\footnotesize (a) Linear $F$: $x$ and $y$ versus $z$
\end{minipage}\hfill
\begin{minipage}{0.48\linewidth}
\centering\includegraphics[width=0.55\linewidth]{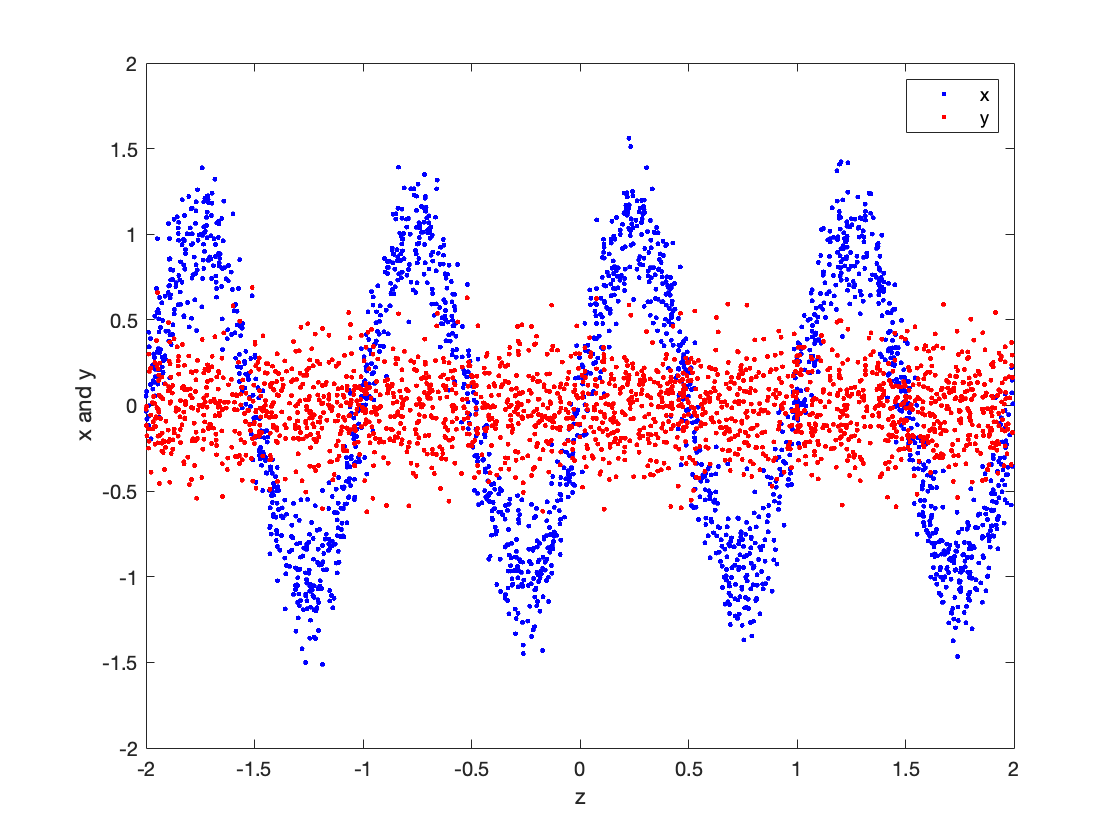}\par
\footnotesize (b) Fixed nonlinear $F$: $x$ and $y$ versus $z$
\end{minipage}

\vspace{0.4em}
\begin{minipage}{0.48\linewidth}
\centering\includegraphics[width=0.55\linewidth]{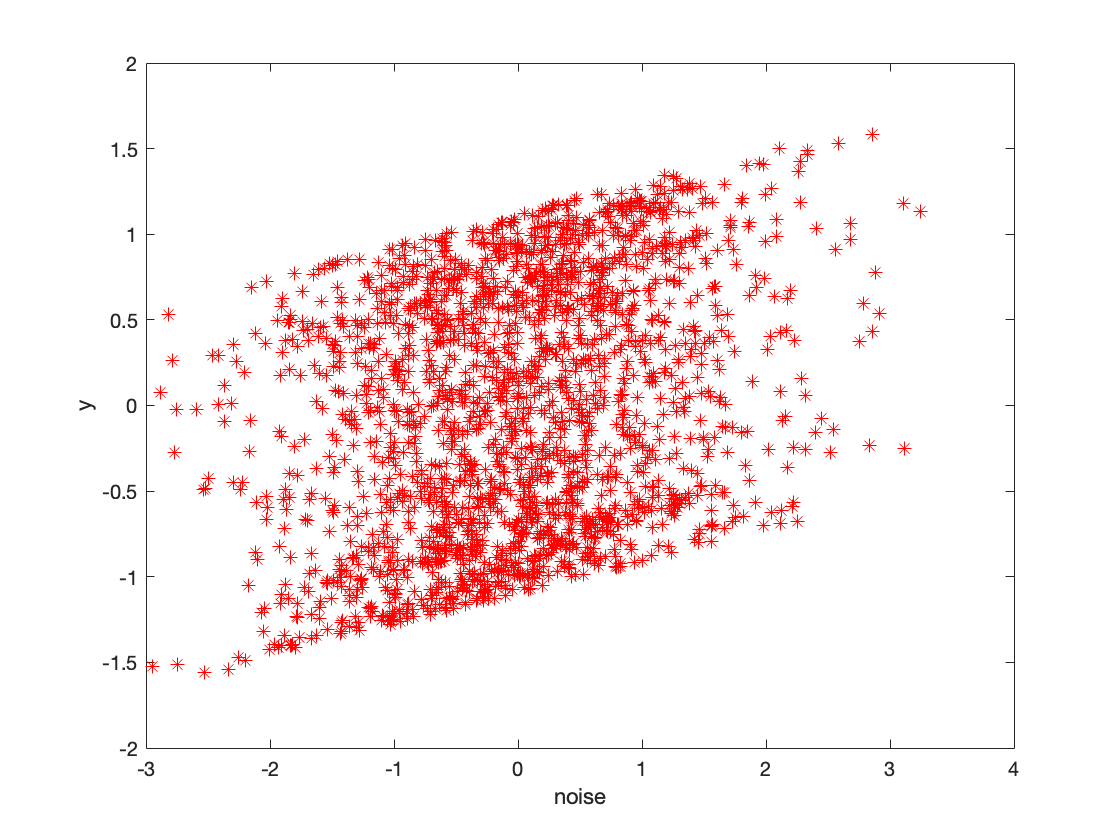}\par
\footnotesize (c) Linear $F$: noise versus $y$
\end{minipage}\hfill
\begin{minipage}{0.48\linewidth}
\centering\includegraphics[width=0.55\linewidth]{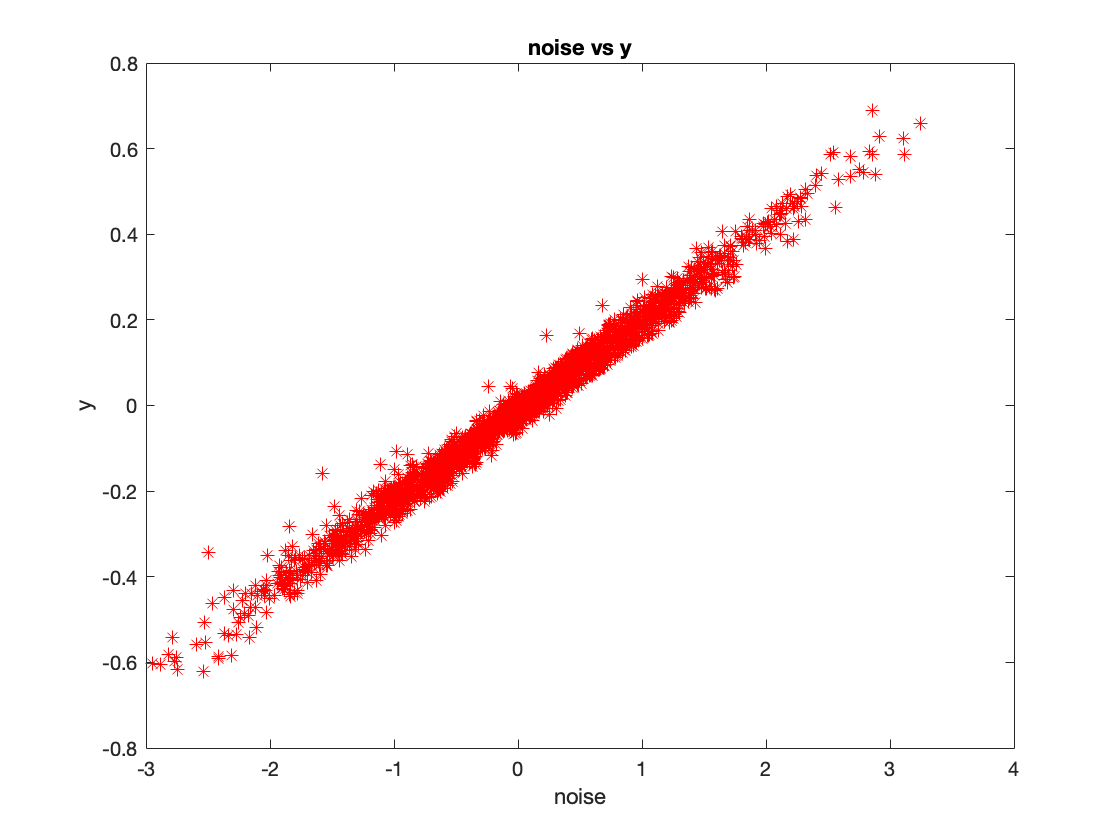}\par
\footnotesize (d) Fixed nonlinear $F$: noise versus $y$
\end{minipage}
\caption{Sinusoidal mean-dependence example. The linear feature matrix cannot
represent the oscillatory dependence, while the fixed nonlinear RBF feature
matrix reduces the residual signal correlation.}
\label{fig:sinusoidal-mean}
\end{figure}

\paragraph{Sinusoidal mean-dependence example.}
We next consider
\begin{equation}
x=\sin(2\pi z)+0.2\epsilon,\qquad
\epsilon\sim\mathcal{N}(0,1).
\label{eq:sinusoidal-mean}
\end{equation}
This example isolates nonlinear mean dependence. The visual comparison in
Figure~\ref{fig:sinusoidal-mean} confirms that the nonlinear feature matrix
removes the sinusoidal structure that is invisible to the linear features.

\begin{table}[ht]
\centering
\caption{Residual dependence after the OTBP update using linear and fixed
nonlinear RBF feature matrices.}
\begin{tabular}{lcc}
\hline
Example and metric & Linear $F$ & Fixed nonlinear $F$ \\
\hline
Nonlinear heteroscedastic:
$\operatorname{corr}(y,z)$
& $0.00414$ & $0.01128$ \\

Nonlinear heteroscedastic:
$\operatorname{corr}((y-\bar y)^2,\sigma(z)^2)$
& $0.41926$ & $0.01575$ \\

Sinusoidal:
$\operatorname{corr}(y,\sin(2\pi z))$
& $0.93421$ & $0.00293$ \\
\hline
\end{tabular}
\label{tab:fixed_nonlinear}
\end{table}

Together, the two examples show why the nonlinear feature matrix is useful:
it detects nonlinear changes in both the conditional variance and the
conditional mean. 

The next section keeps the same nonlinear representation but
allows its centers and bandwidths to be updated, so that comparable dependence
removal can be achieved with fewer centers.

\section{Optimization of RBF Feature Parameters}
\label{sec:alternating-nonlinear}

\subsection{Nonlinear Basis updates}

In the fixed nonlinear setting, the RBF centers and bandwidths are chosen before the OTBP iterations and remain unchanged. Although Gaussian RBFs provide a flexible representation of nonlinear functions, the effectiveness of a finite RBF dictionary depends on the placement of its centers and the choice of bandwidths. We therefore allow these parameters to adapt by optimizing the Step~1 dependence objective to improve the detection of dependence within a limited feature budget.

For the per-center RBF representation,
\begin{equation}
F^{\mathrm{raw}}_{ij}
=
\exp\left[
-\frac{(z_i-z_j^c)^2}{2h_j^2}
\right],
\qquad
j=1,\ldots,n_c,
\end{equation}
where $z_j^c$ and $h_j$ denote the center and bandwidth of the $j$th RBF
feature.

Recall that the Step~1 objective is
\[
J(f,g)
=
\frac{1}{n}f^\top g
-
\frac{1}{2n^2}\|f\|^2\|g\|^2.
\]
Its derivative with respect to $f$ is
\begin{equation}
\frac{\partial J}{\partial f}
=
\frac{1}{n}g
-
\frac{\|g\|^2}{n^2}f.
\end{equation}

\paragraph{Center update.}
For the raw RBF feature,
\begin{equation}
\frac{\partial F^{\mathrm{raw}}_{ij}}
{\partial z_j^c}
=
F^{\mathrm{raw}}_{ij}
\frac{z_i-z_j^c}{h_j^2}.
\end{equation}
Because the actual feature matrix is column-centered, the derivative is
centered in the same way:
\begin{equation}
\frac{\partial F_{ij}}{\partial z_j^c}
=
\frac{\partial F^{\mathrm{raw}}_{ij}}{\partial z_j^c}
-
\frac{1}{n}
\sum_{r=1}^n
\frac{\partial F^{\mathrm{raw}}_{rj}}{\partial z_j^c}.
\end{equation}
Since
\[
f_i=\sum_{j=1}^{n_c}a_jF_{ij},
\]
we have
\[
\frac{\partial f_i}{\partial z_j^c}
=
a_j
\frac{\partial F_{ij}}{\partial z_j^c}.
\]
Therefore,
\begin{equation}
\frac{\partial J}{\partial z_j^c}
=
\sum_{i=1}^n
\left(
\frac{1}{n}g_i
-
\frac{\|g\|^2}{n^2}f_i
\right)
a_j
\frac{\partial F_{ij}}{\partial z_j^c}.
\end{equation}

Since Step~1 maximizes $J$, the centers are updated by gradient ascent:
\begin{equation}
z_j^c
\leftarrow
z_j^c
+
\eta_c
\frac{\partial J}{\partial z_j^c}.
\end{equation}
After each accepted update, the feature matrix $F$ is reconstructed and
re-centered, and the coefficient vectors are recomputed.

\paragraph{Per-center bandwidth update.}
Instead of using one common bandwidth for all RBF features, we allow each
center to have its own bandwidth,
\[
h=(h_1,\ldots,h_{n_c})^\top.
\]
For each $h_j$,
\begin{equation}
\frac{\partial F^{\mathrm{raw}}_{ij}}
{\partial h_j}
=
F^{\mathrm{raw}}_{ij}
\frac{(z_i-z_j^c)^2}{h_j^3}.
\end{equation}
After centering the derivative in the same way as $F$, the chain rule gives
\begin{equation}
\frac{\partial J}{\partial h_j}
=
\sum_{i=1}^n
\frac{\partial J}{\partial f_i}
a_j
\frac{\partial F_{ij}}{\partial h_j}.
\end{equation}

Collecting the derivatives for all centers gives the bandwidth gradient
$\nabla_h J$. The bandwidth vector is then updated by
\begin{equation}
h
\leftarrow
h
+
\eta_h \nabla_h J.
\end{equation}

\paragraph{Armijo-type accept--reject rule.}
Both the center and bandwidth updates use the same Armijo-type
accept--reject strategy. Let $\theta$ denote the current block of basis
parameters, either the center vector $z^c$ or the bandwidth vector $h$.
Given the current gradient $\nabla_\theta J$, we propose the gradient-ascent
trial update
\[
\theta_{\mathrm{trial}}
=
\theta_{\mathrm{old}}
+
\eta_\theta \nabla_\theta J.
\]

Using the trial parameters, we reconstruct and center the feature matrix
$F$, recompute the Step~1 coefficient vectors, and evaluate
\[
J_{\mathrm{new}}
=
J(f_{\mathrm{trial}},g_{\mathrm{trial}}).
\]
The trial update is accepted if
\[
J_{\mathrm{new}}-J_{\mathrm{old}}
>
c_\theta \eta_\theta
\|\nabla_\theta J\|^2
-
\varepsilon,
\]
where $c_\theta>0$ is the Armijo parameter and $\varepsilon$ is a small numerical tolerance.

If the condition is satisfied, the trial update is accepted and the
learning rate is slightly increased:
\[
\theta
\leftarrow
\theta_{\mathrm{trial}},
\qquad
\eta_\theta
\leftarrow
1.01\,\eta_\theta.
\]
Otherwise, the trial update is rejected and the learning rate is reduced:
\[
\theta
\leftarrow
\theta_{\mathrm{old}},
\qquad
\eta_\theta
\leftarrow
0.81\,\eta_\theta.
\]

\paragraph{Alternating center and bandwidth updates.}
The center and bandwidth updates described above can be combined in a
block-wise alternating scheme. Starting from the current RBF basis, we first
update the centers $z^c$ while keeping the bandwidths $h$ fixed, and then
recompute and center the feature matrix $F$. Next, with the updated centers
held fixed, we update the bandwidth vector $h$ and reconstruct $F$ again.
The procedure is summarized as

\[
\begin{array}{l}
\text{initialize } z^c \text{ and } h, \\[2mm]
\text{for each Step~1 dependence-search step:} \\[1mm]
\qquad \text{for } b=1,\ldots,n_{\mathrm{alt}}: \\[1mm]
\qquad\qquad \text{update } z^c \text{ with } h \text{ fixed},\\
\qquad\qquad \text{recompute and center } F,\\
\qquad\qquad \text{update } h \text{ with } z^c \text{ fixed},\\
\qquad\qquad \text{recompute and center } F.
\end{array}
\]

Both blocks maximize the same Step~1 objective $J$ and use an Armijo-type
accept--reject rule. The alternating updates are repeated for a prescribed
number of blocks or until the parameter movement and objective improvement
become sufficiently small. After the basis update is completed, the
coefficient vectors are recomputed and the algorithm proceeds to Step~2.

\subsection{Numerical evaluation of basis updates}

We compare a fixed RBF basis with the alternating basis update on the
sinusoidal example. 

\begin{table}[ht]
\centering
\caption{Comparison of fixed and alternating RBF bases for the sinusoidal
example.}
\begin{tabular}{lccc}
\hline
Method & $n_c$ & $\operatorname{corr}(y,\epsilon)$
& $\operatorname{corr}(y,\sin(2\pi z))$ \\
\hline
Fixed centers and bandwidth
& 10 & 0.64755 & 0.29816 \\

Alternating centers and bandwidths
& 10 & 0.96414 & 0.038134 \\

Fixed centers and bandwidth
& 20 & 0.99260 & 0.020439 \\
\hline
\end{tabular}
\label{tab:fixed_vs_alternating}
\end{table}

With the same number of RBF features, the alternating update substantially
improves the result: for $n_c=10$, the residual sinusoidal correlation
decreases from $0.29816$ for the fixed basis to $0.038134$ after adapting
the centers and bandwidths. 

These results indicate that basis updates can compensate for a relatively
small feature space. In this example, the alternating update with
$n_c=10$ reaches a low-residual regime comparable to that of a fixed basis
with $n_c=20$.

\begin{figure}[H]
\centering

\begin{minipage}{0.40\linewidth}
\centering\includegraphics[width=0.74\linewidth]{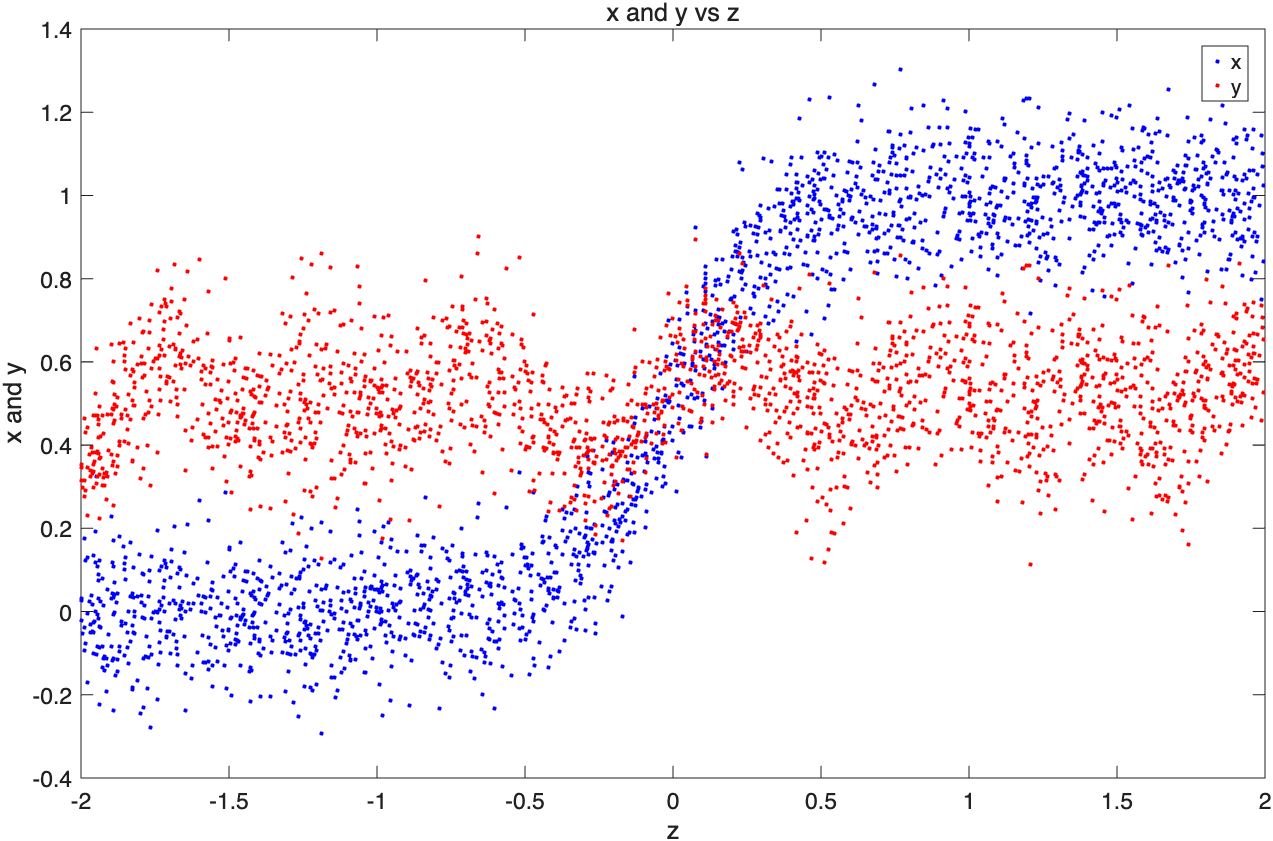}\par
\footnotesize (a) Fixed: $x$ and $y$ versus $z$
\end{minipage}
\hspace{0.05\linewidth}
\begin{minipage}{0.40\linewidth}
\centering\includegraphics[width=0.74\linewidth]{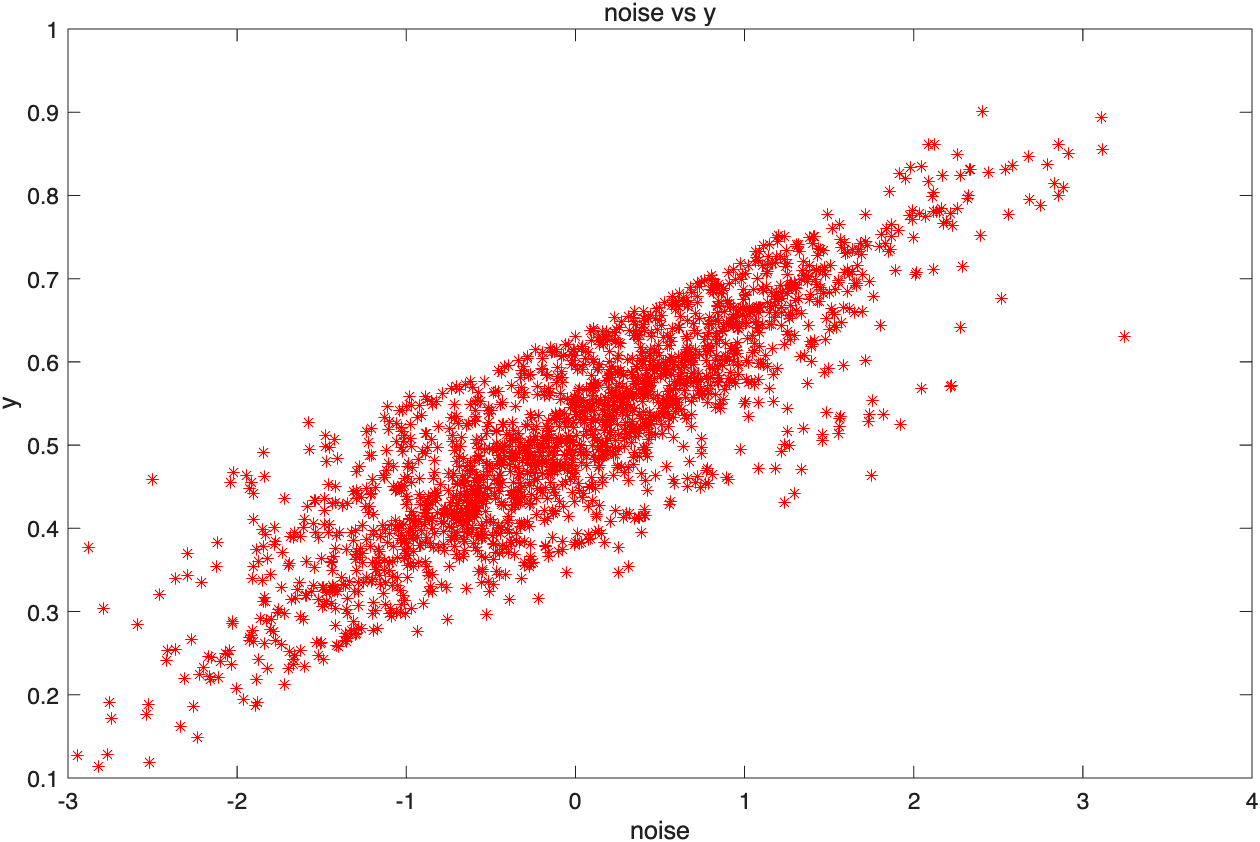}\par
\footnotesize (b) Fixed: noise versus $y$
\end{minipage}

\vspace{0.7em}

\begin{minipage}{0.31\linewidth}
\centering\includegraphics[width=0.74\linewidth]{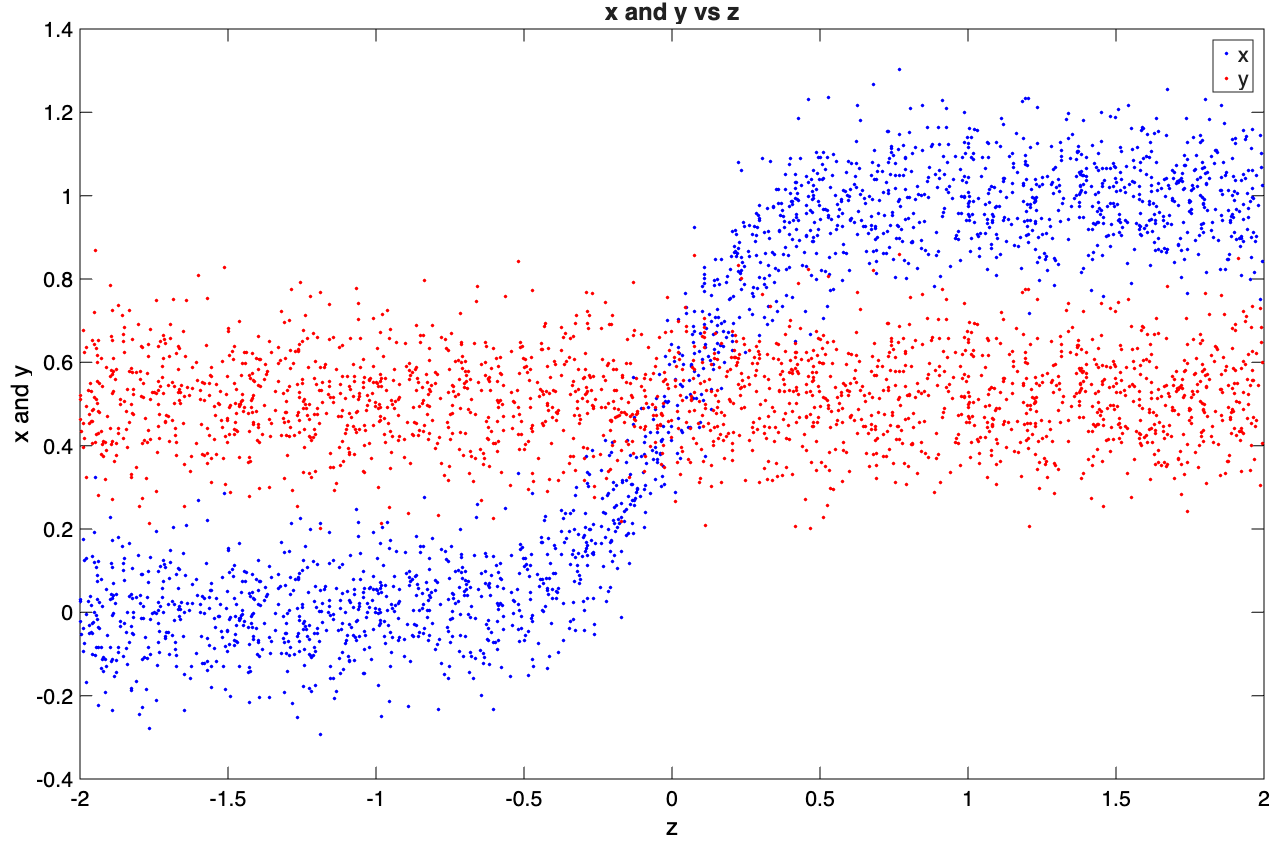}\par
\footnotesize (c) Alternating: $x$ and $y$ versus $z$
\end{minipage}
\hfill
\begin{minipage}{0.31\linewidth}
\centering\includegraphics[width=0.74\linewidth]{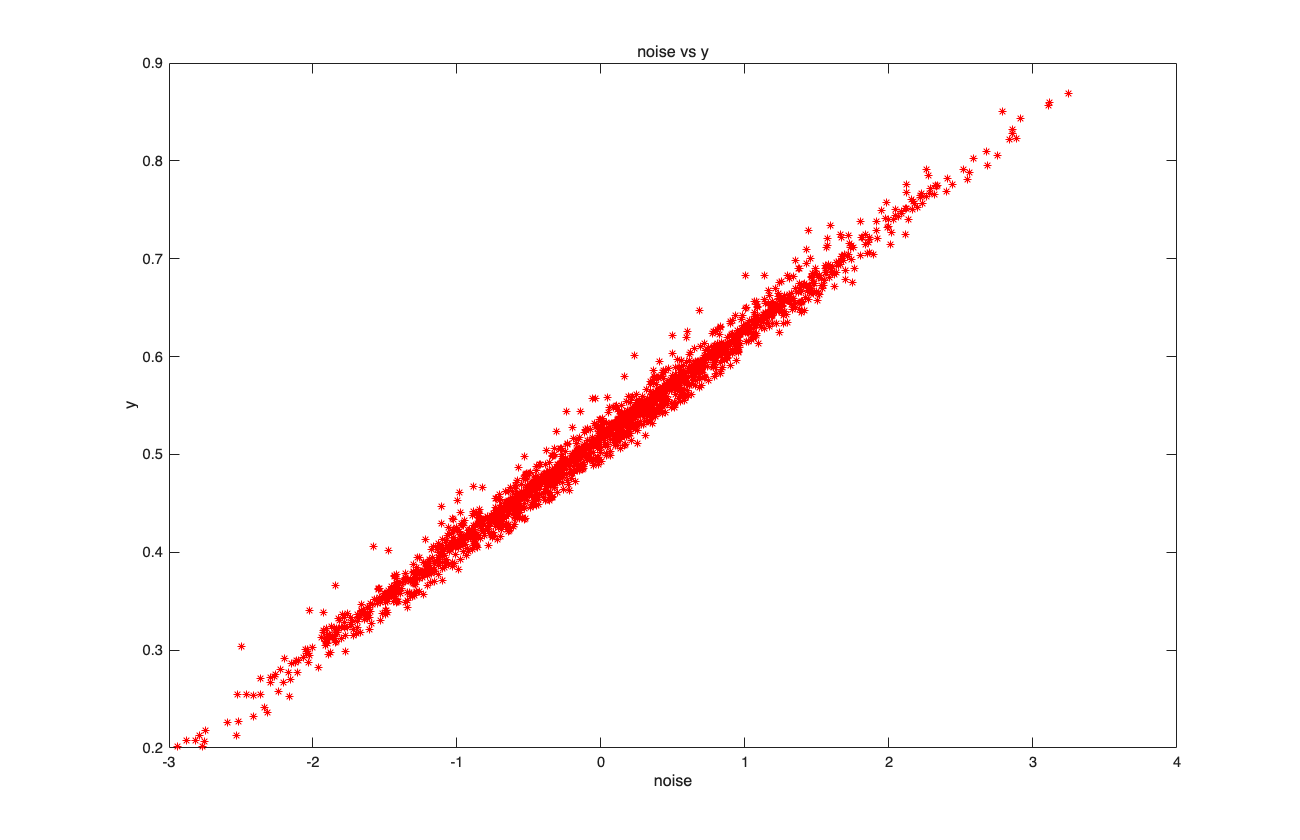}\par
\footnotesize (d) Alternating: noise versus $y$
\end{minipage}
\hfill
\begin{minipage}{0.31\linewidth}
\centering\includegraphics[width=0.74\linewidth]{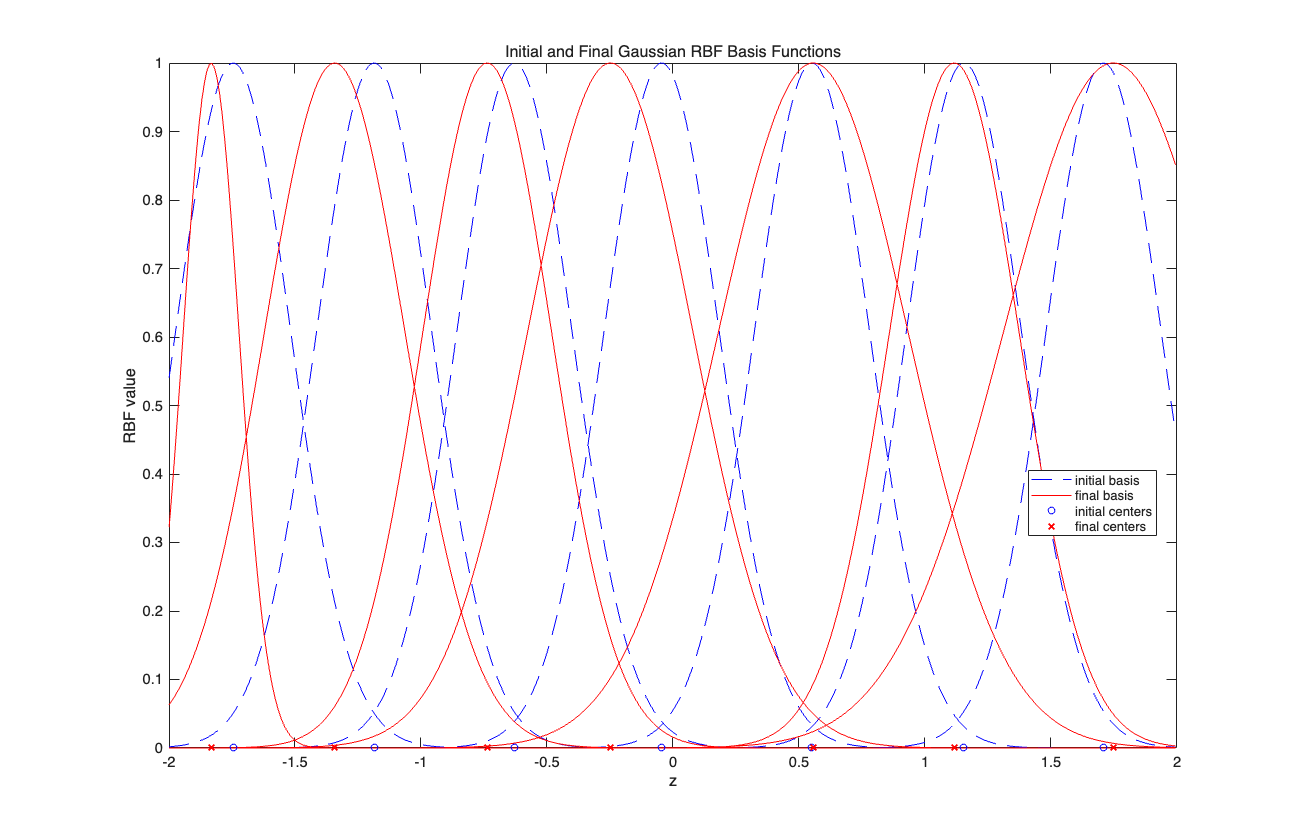}\par
\footnotesize (e) Initial and final RBF basis
\end{minipage}

\caption{Sigmoid example with $n_c=7$. The top row shows the result obtained
with the fixed RBF basis. The bottom row shows the result after alternating
center and bandwidth updates, together with the corresponding evolution of
the RBF basis.}
\label{fig:sigmoid-nc7-comparison}
\end{figure}

Next, we consider more examples.
\paragraph{Sigmoid example.}
We next consider the monotone nonlinear dependence
\[
x=\frac{1}{1+\exp(-6z)}+\epsilon.
\]
With only $n_c=7$ centers, the fixed basis leaves residual sigmoid correlation
$0.092653$, whereas the alternating center-bandwidth update reduces it to
$0.021301$.  At the same time, the correlation between $y$ and the noise
increases from $0.85694$ to $0.99391$.  The visual comparison in
Figure~\ref{fig:sigmoid-nc7-comparison} shows that adapting the same seven
features removes the nonlinear transition more completely.

\begin{figure}[H]
\centering

\begin{minipage}{0.40\linewidth}
\centering\includegraphics[width=0.74\linewidth]{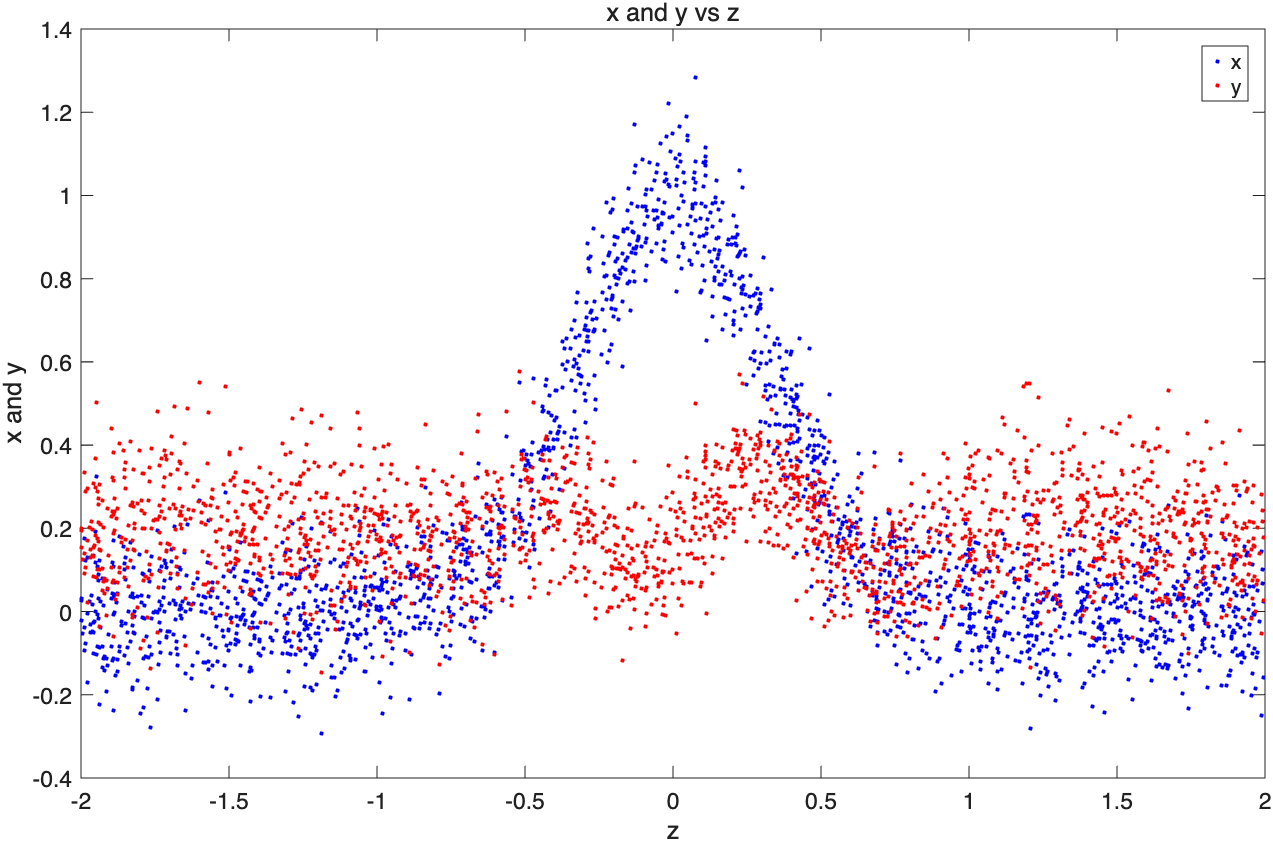}\par
\footnotesize (a) Fixed: $x$ and $y$ versus $z$
\end{minipage}
\hspace{0.05\linewidth}
\begin{minipage}{0.40\linewidth}
\centering\includegraphics[width=0.74\linewidth]{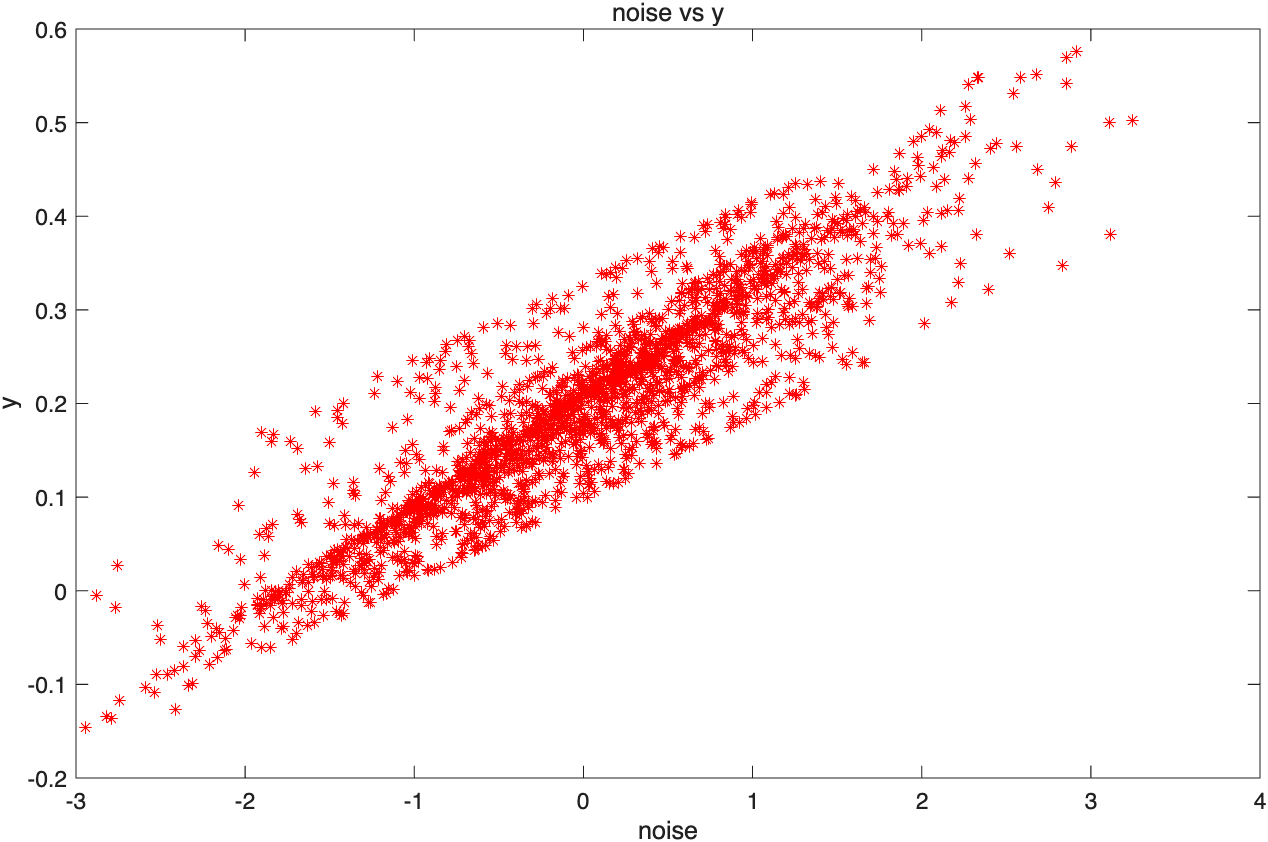}\par
\footnotesize (b) Fixed: noise versus $y$
\end{minipage}

\vspace{0.7em}

\begin{minipage}{0.31\linewidth}
\centering\includegraphics[width=0.74\linewidth]{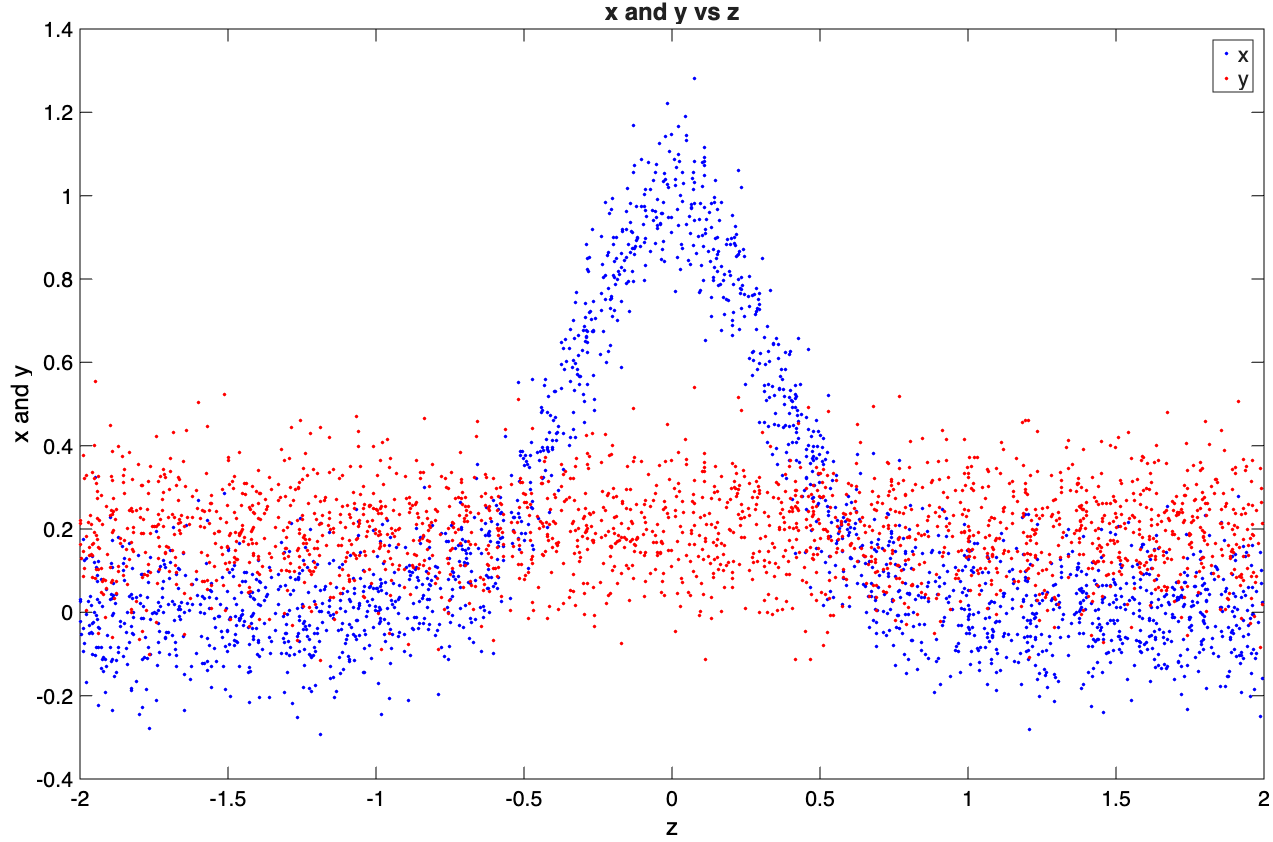}\par
\footnotesize (c) Alternating: $x$ and $y$ versus $z$
\end{minipage}
\hfill
\begin{minipage}{0.31\linewidth}
\centering\includegraphics[width=0.74\linewidth]{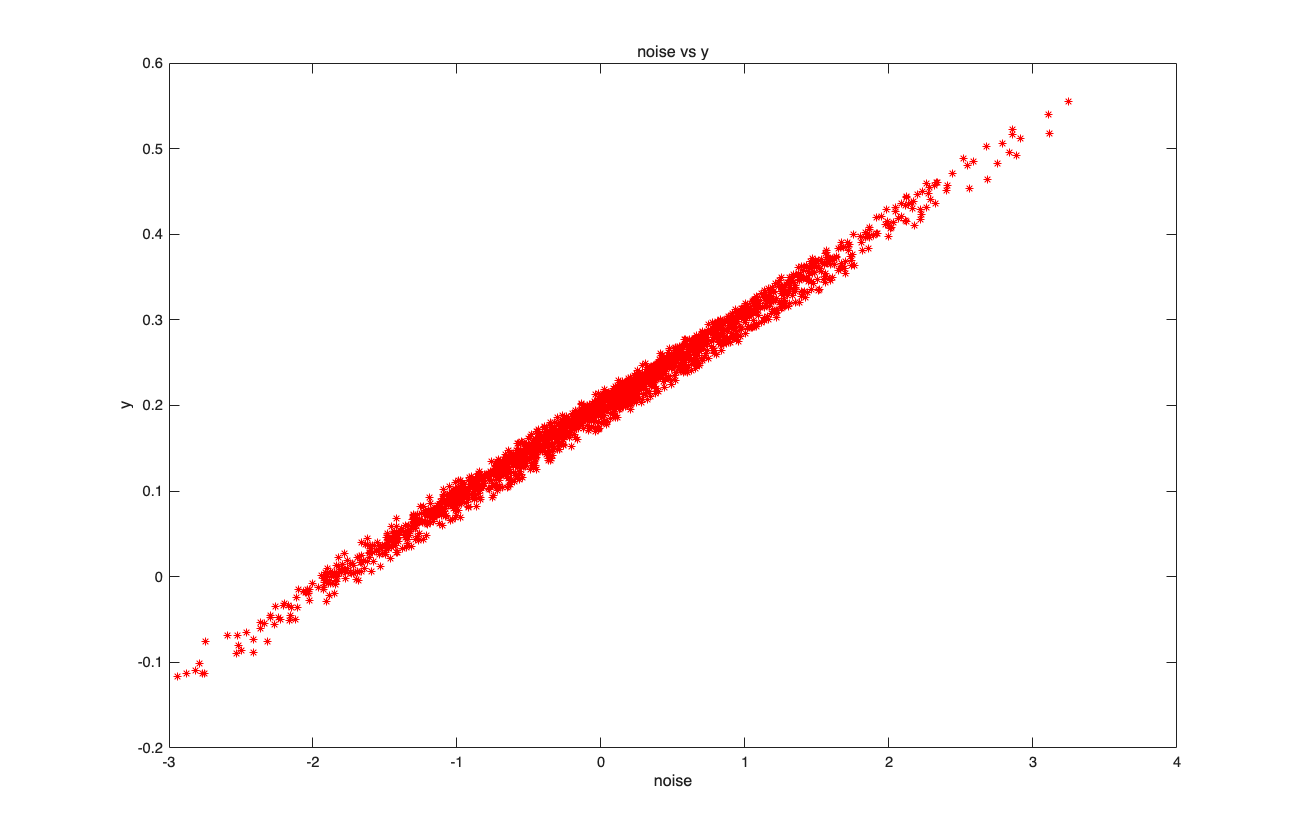}\par
\footnotesize (d) Alternating: noise versus $y$
\end{minipage}
\hfill
\begin{minipage}{0.31\linewidth}
\centering\includegraphics[width=0.74\linewidth]{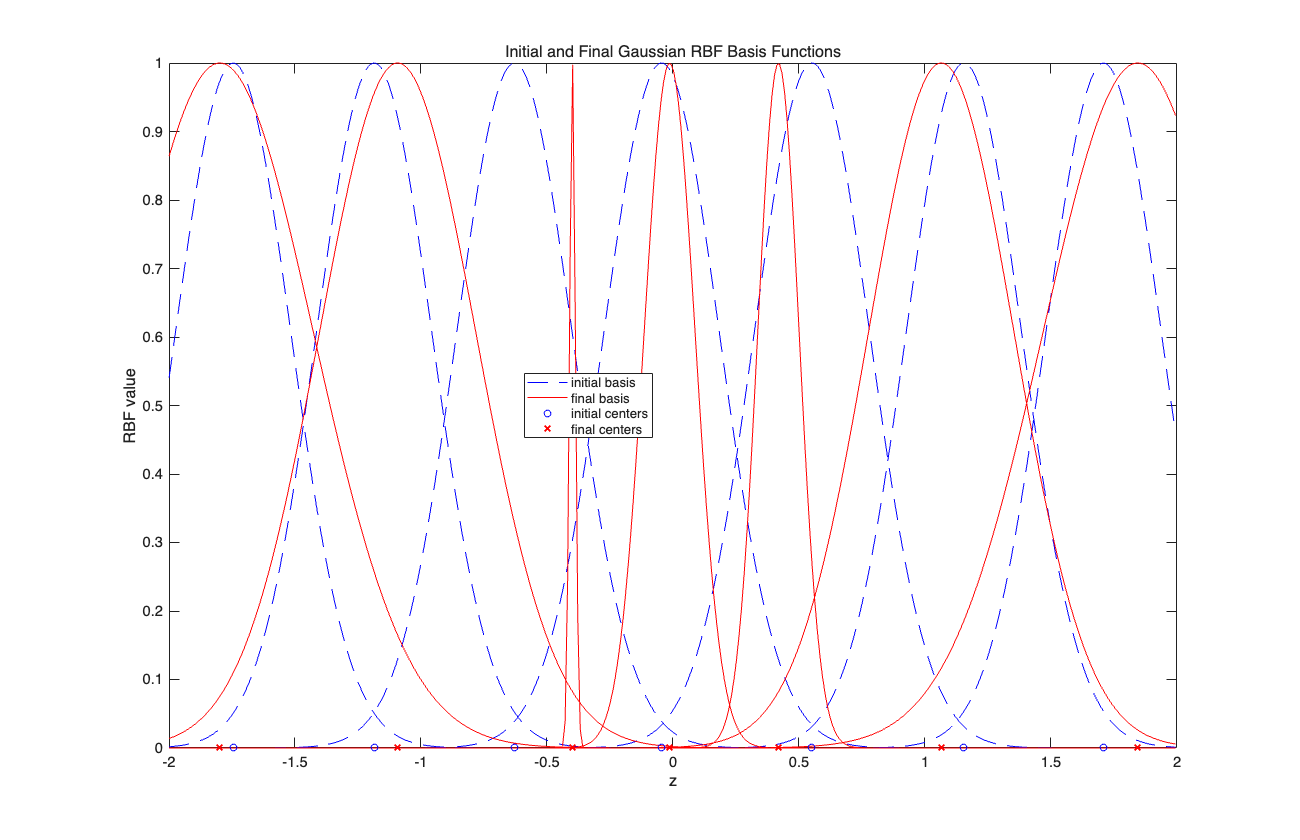}\par
\footnotesize (e) Initial and final RBF basis
\end{minipage}

\caption{Localized middle-bump example with $n_c=7$. The top row shows the
fixed-basis result. The bottom row shows the transformed samples after
alternating center and bandwidth updates and the corresponding evolution of
the RBF basis.}
\label{fig:middlebump-nc7-comparison}
\end{figure}

\paragraph{Localized middle-bump example.}
Finally, we consider
\[
x=\exp\!\left(-\frac{z^2}{0.2}\right)+0.1\epsilon,
\]
whose dependence is concentrated near the center of the domain.  Again using
only $n_c=7$, the fixed basis gives residual middle-bump correlation
$0.086089$, while the alternating update reduces it to $0.020615$.  The noise
correlation increases from $0.90597$ to $0.99542$.  As shown in
Figure~\ref{fig:middlebump-nc7-comparison}, the alternating basis is better
able to represent and remove the localized dependence without increasing the
number of centers.

\begin{table}[H]
\centering
\small
\caption{Fixed and alternating RBF bases for the sigmoid and middle-bump
examples.  All comparisons use $n_c=7$.}
\label{tab:nc7-additional-examples}
\begin{tabular}{llcc}
\toprule
Example & Method & $\operatorname{corr}(y,\epsilon)$ & Residual correlation\\
\midrule
Sigmoid & Fixed & 0.85694 & 0.092653\\
Sigmoid & Alternating & 0.99391 & 0.021301\\
Middle bump & Fixed & 0.90597 & 0.086089\\
Middle bump & Alternating & 0.99542 & 0.020615\\
\bottomrule
\end{tabular}
\end{table}

The sigmoid and middle-bump experiments show that the benefit of adaptation is
not limited to periodic signals: it also appears for monotone and localized
nonlinear structure.

\subsection{Controlled validation of the basis updates}
\label{sec:controlled-basis-validation}

The preceding examples show that adapting the nonlinear feature map can improve the transformed samples.  We next isolate Step~1 in order to examine more directly whether the center and bandwidth updates optimize the intended dependence criterion.  Throughout this subsection, Step~2 is frozen by setting $y=x$.  Consequently, the dependence pattern does not change during the optimization, and the behavior of the basis update can be compared with a known target.

\subsubsection{Single-center ground-truth recovery}

We consider the single-bump model
\begin{equation}
    z_i\sim\operatorname{Unif}[-2,2],\qquad
    s_i=\exp\!\left[-\frac{(z_i-c^\star)^2}{2(h^\star)^2}\right],\qquad
    x_i=s_i+0.1\varepsilon_i,
    \label{eq:known-bump}
\end{equation}
where $\varepsilon_i\sim\mathcal N(0,1)$ and $n=3000$.  Since the signal is itself a Gaussian radial basis function, its parameters are known exactly.  In the centered case,
\begin{equation}
    c^\star=0,
    \qquad
    h^\star=\sqrt{0.1}\simeq0.3162.
\end{equation}
For every trial value of the basis parameters, the coefficients $a$ and $b$ are re-optimized by the same \texttt{fg} routine used in the full algorithm.  We write the resulting profiled criterion as
\begin{equation}
    J^\star(c,h)=\max_{a,b}J\bigl(F(z;c,h)a,G(y)b\bigr).
    \label{eq:profiled-J}
\end{equation}
The parameter updates use the original center and bandwidth gradients and the same acceptance rule as the full OTBP implementation.

\begin{figure}[H]
    \centering
    \includegraphics[width=0.92\textwidth]{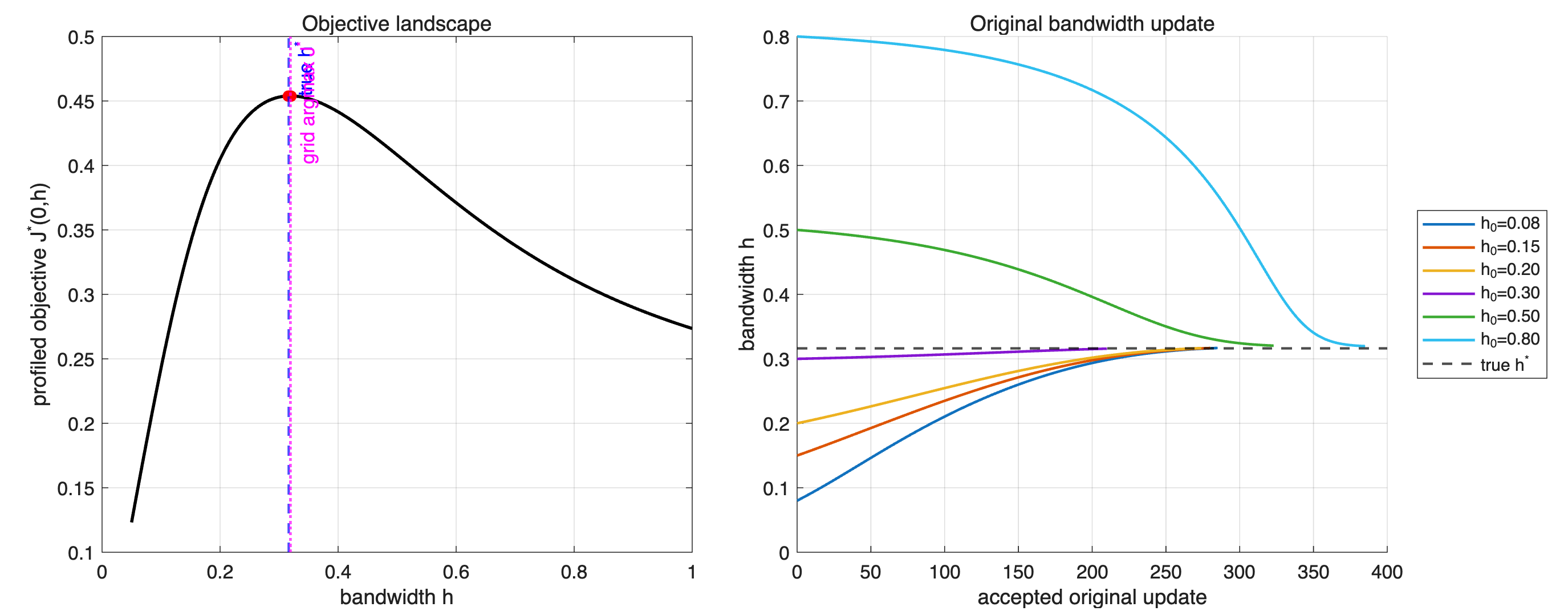}
    \caption{Single-center bandwidth recovery with Step~2 frozen.  The left panel shows the profiled objective $J^\star(0,h)$; the right panel shows the original bandwidth-update trajectories from six initial values.  The finite-sample maximizer is close to the generative bandwidth $h^\star=\sqrt{0.1}$.}
    \label{fig:single-bandwidth-recovery}
\end{figure}

\paragraph{Bandwidth recovery.}
We first fix the center at $c^\star$ and initialize the bandwidth at
\[
    h_0\in\{0.08,0.15,0.20,0.30,0.50,0.80\}.
\]
All six trajectories converge to the interval $[0.3158,0.3205]$, close to both the generative value $h^\star=0.3162$ and the finite-sample grid maximizer of $J^\star(0,h)$.  The final objective values lie between $0.45393$ and $0.45396$.  Figure~\ref{fig:single-bandwidth-recovery} shows that the update approaches the same maximizer from both narrow and broad initial bandwidths.

\paragraph{Center recovery.}
We then fix $h=h^\star$ and update only the center.  Initializations at $c_0=-0.5$ and $c_0=0.3$ converge to $c=0.0057$ and $c=0.0069$, respectively, and attain $J\simeq0.4541$.  However, initializations farther from the signal converge to secondary stationary points near $c=-1.37$ and $c=1.40$, where the objective is substantially smaller.  Thus, the center update follows the correct local ascent direction, but the center landscape is nonconvex and its basin of attraction does not cover the full domain.

\begin{figure}[H]
    \centering
    \includegraphics[width=0.92\textwidth]{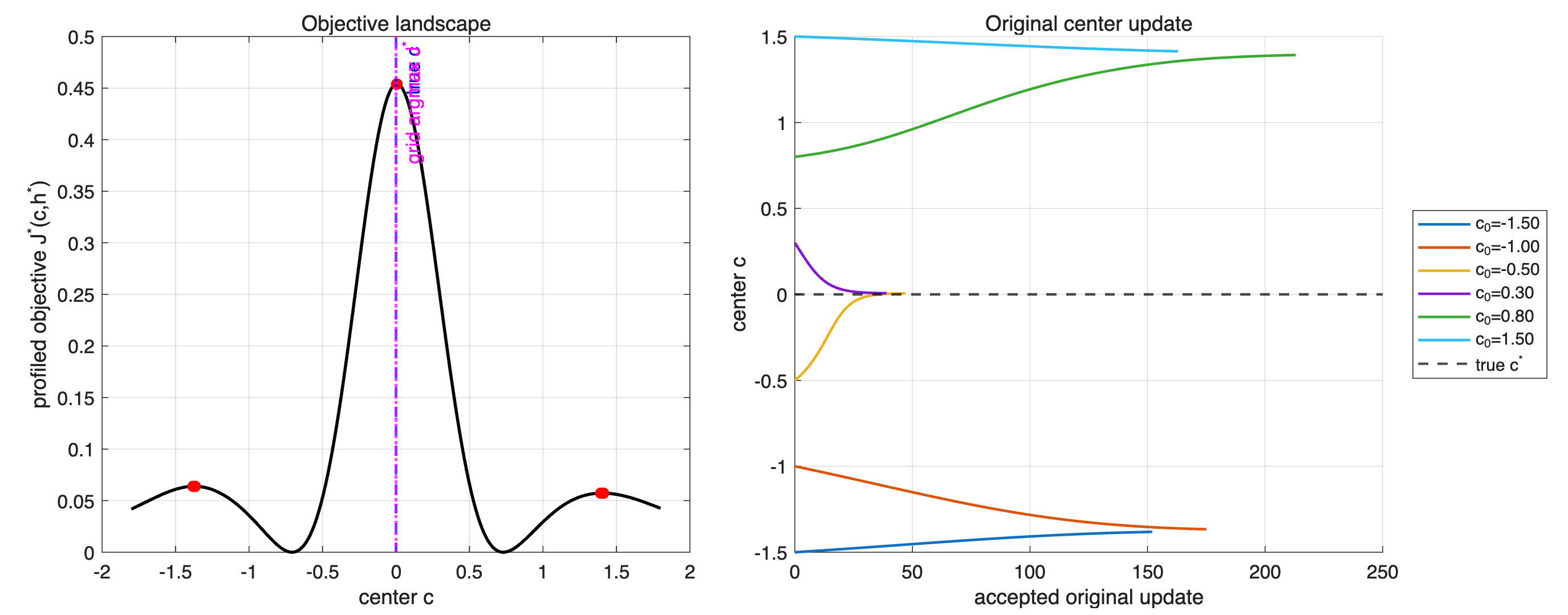}
    \caption{Single-center recovery with the bandwidth fixed at $h^\star$.  Starts within the central basin converge to the true center, whereas distant starts converge to lower local maxima of the profiled objective.}
    \label{fig:single-center-recovery}
\end{figure}

\begin{figure}[H]
    \centering
    \includegraphics[width=0.50\textwidth]{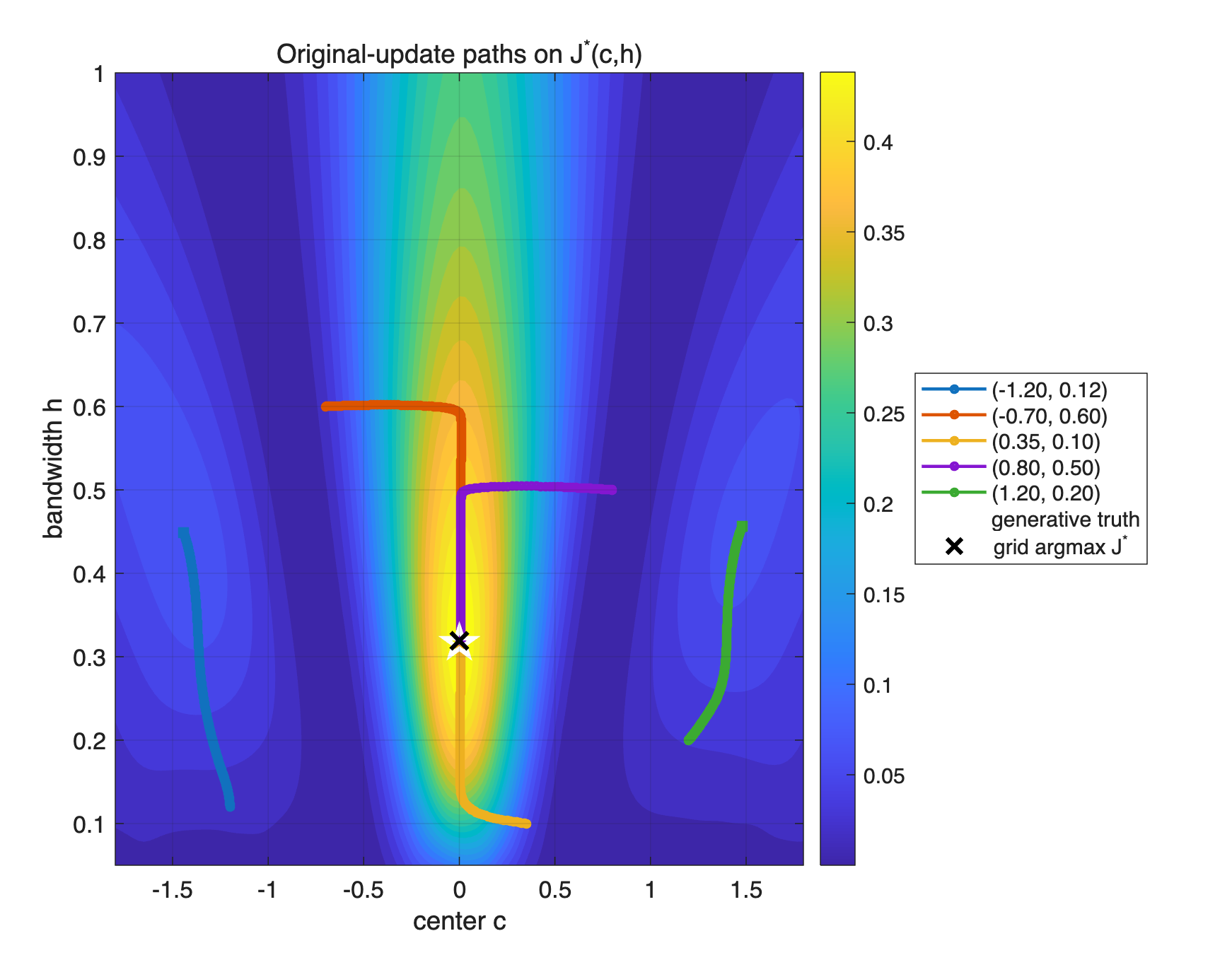}
    \caption{Alternating center--bandwidth trajectories superimposed on the profiled objective $J^\star(c,h)$.  The black cross denotes the finite-sample grid maximizer and the star denotes the generative parameters.  The side basins explain why a single distant initialization can fail.}
    \label{fig:single-joint-landscape}
\end{figure}

\paragraph{Joint alternating recovery.}
Finally, we update the center and bandwidth alternately.  The joint trajectories in Figure~\ref{fig:single-joint-landscape} confirm the same structure.  Starts in the central basin converge close to $(c^\star,h^\star)$ and attain $J\simeq0.4541$, while distant starts remain in low-objective side basins.  This controlled experiment therefore verifies the update equations near the desired solution while also identifying initialization sensitivity as the main limitation of the one-center problem.

\subsubsection{Broad-bandwidth initialization.}
The preceding experiments reveal an initialization limitation of the
single-center problem. With a relatively narrow RBF, center
initializations sufficiently far from the true dependence region can
converge to low-objective local maxima instead of the global solution at
$c^\star=0$. This motivates the use of a broader initial bandwidth,
which allows a single basis function to detect dependence over a wider
region of the domain.

\begin{figure}[t]
    \centering
    \includegraphics[width=0.85\textwidth]{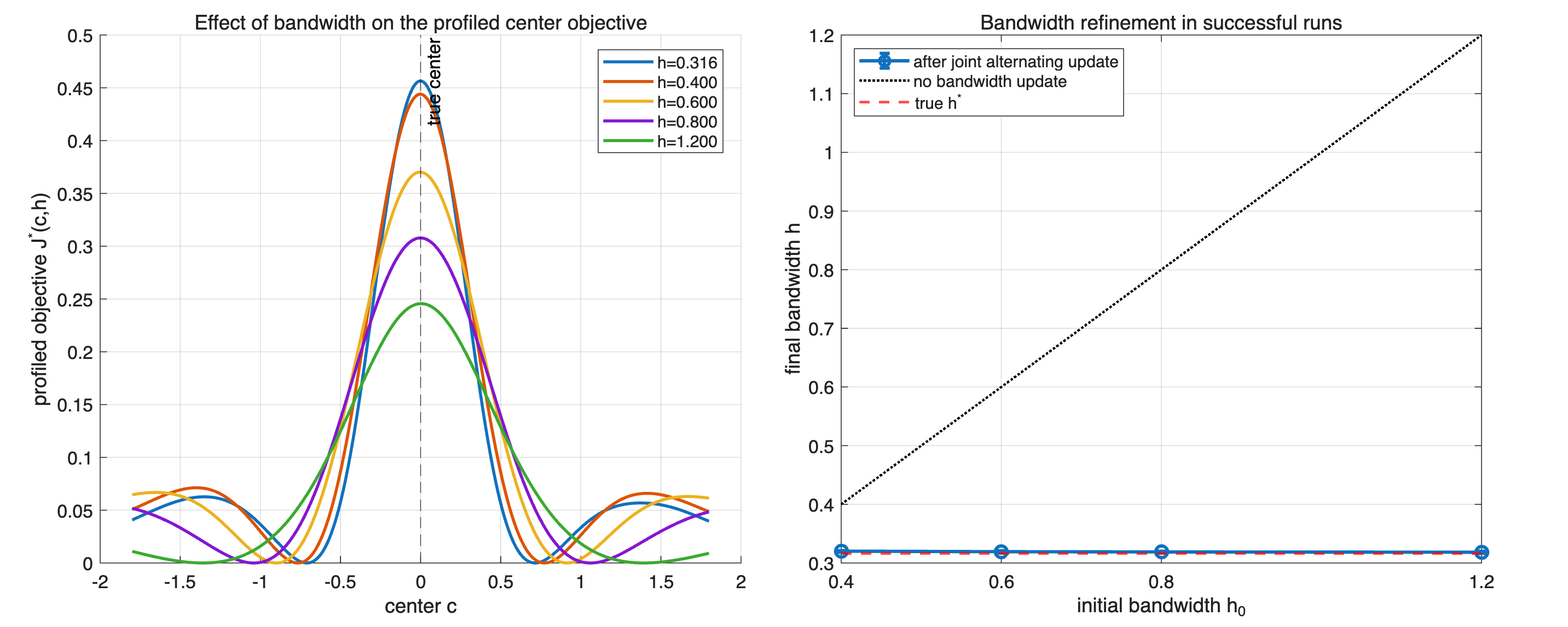}
    \caption{Effect of broad-bandwidth initialization in the
    single-center ground-truth experiment. Left: increasing the fixed
    bandwidth broadens and smooths the central attraction basin of the
    profiled objective, allowing center updates from a wider range of
    initializations to reach the true center. However, sufficiently
    distant starts can still converge to side maxima. Right: among the
    successful runs, the joint alternating updates refine each broad
    initial bandwidth toward the true value
    $h^\star=\sqrt{0.1}\approx0.3162$. The points and error bars summarize
    the results over ten independently generated data sets.}
    \label{fig:broad-bandwidth}
\end{figure}

Figure~\ref{fig:broad-bandwidth} illustrates the resulting mechanism.
The left panel shows the profiled center objective
\[
J^\star(c,h)=\max_{a,b}J\bigl(F(z;c,h)a,G(y)b\bigr)
\]
for several fixed bandwidths. Increasing $h$ broadens and smooths the
central attraction basin around the true center. With the original
narrow initialization, only starts close to the dependence region, such
as $|c_0|\leq 0.5$, reliably reach the central solution. With
$h_0=0.8$ or $1.2$, the recoverable initialization range expands to
$|c_0|=1$. The more distant starts $c_0=\pm1.5$ remain outside the
central basin, showing that a broad bandwidth improves robustness but
does not make the objective globally convex.

After locating the central dependence region, the bandwidth should not
remain artificially broad. We therefore perform a second, independent
experiment in which the center and bandwidth are updated alternately
from the same initial pairs $(c_0,h_0)$. The right panel of
Figure~\ref{fig:broad-bandwidth} summarizes the successful runs over ten
independently generated data sets. Although the initial bandwidth ranges
from $0.4$ to $1.2$, the alternating updates consistently refine it to
approximately
\[
h_{\mathrm{final}}\approx0.318\text{--}0.320,
\]
which is close to the generative value
\[
h^\star=\sqrt{0.1}\approx0.3162.
\]
Thus, the broad initialization and the bandwidth update have
complementary roles: the former enlarges the attraction basin for the
center search, whereas the latter recovers the localized spatial scale
after the correct dependence region has been identified.

\newpage

\subsection{Robustness to Initialization and Covariate Distribution}

The previous experiments use uniformly distributed covariates, for which the
initial placement of the RBF centers provides relatively even coverage of the
domain. To examine whether the behavior changes when the distribution of
$z$ is non-uniform, we compare two covariate distributions:
\[
z\sim \mathrm{Uniform}[-2,2],
\]
and a non-uniform mixture-Gaussian distribution. For each distribution, we
consider three different center initializations: $k$-means, a uniform grid,
and random initialization. In all cases, we use $n_c=7$ RBF centers and then
apply the same alternating center--bandwidth update described above.

\paragraph{Localized middle-bump dependence.}
We first consider the localized middle-bump example
\[
x
=
\exp\left(-\frac{z^2}{0.2}\right)
+
0.1\epsilon.
\]

\begin{table}[ht]
\centering
\small
\setlength{\tabcolsep}{4pt}
\caption{Effect of $z$ distribution and center initialization for the
middle-bump example with $n_c=7$.}
\begin{tabular}{llcccc}
\hline
$z$ dist. & Init.
& $\operatorname{corr}(y,\mathrm{bump})$
& $\operatorname{corr}(y,\epsilon)$
& $\Delta z^c$
& $\Delta h$ \\
\hline
Uniform
& $k$-means
& 0.031400
& 0.98157
& 0.22292
& 0.0010819 \\

Uniform
& Uniform
& 0.026505
& 0.98910
& 0.26833
& 0.0011222 \\

Uniform
& Random
& 0.026115
& 0.98795
& 0.27726
& 0.0010666 \\

Mixture Gaussian
& $k$-means
& 0.071957
& 0.80388
& 0.65186
& 0.0063813 \\

Mixture Gaussian
& Uniform
& 0.031773
& 0.98731
& 0.57286
& 0.00066421 \\

Mixture Gaussian
& Random
& 0.032272
& 0.98700
& 0.62031
& 0.00069302 \\
\hline
\end{tabular}
\label{tab:middle_bump_initialization}
\end{table}

Here $\Delta z^c$ denotes the total center movement and $\Delta h$ denotes
the total bandwidth movement.

\paragraph{Sigmoid dependence.}
We next consider the sigmoid example
\[
x
=
\frac{1}{1+\exp(-6z)}
+
\epsilon.
\]

\begin{table}[ht]
\centering
\small
\setlength{\tabcolsep}{4pt}
\caption{Effect of $z$ distribution and center initialization for the
sigmoid example with $n_c=7$.}
\begin{tabular}{llcccc}
\hline
$z$ dist. & Init.
& $\operatorname{corr}(y,\mathrm{sigmoid})$
& $\operatorname{corr}(y,\epsilon)$
& $\Delta z^c$
& $\Delta h$ \\
\hline
Uniform
& $k$-means
& 0.030107
& 0.96529
& 0.29171
& $1.3498\times10^{-5}$ \\

Uniform
& Uniform
& 0.032429
& 0.97010
& 0.67419
& $1.2454\times10^{-5}$ \\

Uniform
& Random
& 0.030983
& 0.96900
& 0.35746
& $1.2021\times10^{-5}$ \\

Mixture Gaussian
& $k$-means
& 0.027404
& 0.97962
& 0.89828
& $1.9960\times10^{-5}$ \\

Mixture Gaussian
& Uniform
& 0.027417
& 0.98053
& 0.76299
& 0.0017226 \\

Mixture Gaussian
& Random
& 0.030928
& 0.98036
& 0.75397
& 0.0016786 \\
\hline
\end{tabular}
\label{tab:sigmoid_initialization}
\end{table}

Here $\Delta z^c$ denotes the total center movement and $\Delta h$ denotes
the total bandwidth movement.

For the localized middle-bump example with a mixture-Gaussian distribution of $z$, uniform initialization performs better than $k$-means after the basis updates. This is probably because $k$-means is primarily density-aware: it places more centers in high-density regions, even though lower-density regions may still contain important dependence. In contrast, a uniform initialization provides broader coverage of the $z$-domain, giving the subsequent dependence-driven updates more flexibility to move the centers toward locations that are useful for detecting and removing dependence.

Overall, these results motivate a dependence-aware rather than purely
density-aware basis. By allowing the RBF centers and bandwidths to adapt to the Step~1 dependence objective, the nonlinear feature representation can be useful for a range of nonlinear structures, including oscillatory, monotone, localized, and heteroscedastic dependence. 

\section{Extension to Higher Dimensions}

To explore the extension to higher-dimensional covariates, we first consider the two-dimensional case \(z \in \mathbb{R}^2\). This setting captures the main changes required for the multidimensional formulation while remaining easy to visualize. The extension only changes the geometry of the RBF features and the derivatives with respect to their centers. The Step~1 objective, feature centering, and alternating update procedure remain the same as in the one-dimensional case.

Let
\[
z_i = (z_{i1},z_{i2})^\top,
\qquad
c_j = (c_{j1},c_{j2})^\top \in \mathbb{R}^2.
\]
For each center \(c_j\) with isotropic bandwidth \(h_j>0\), define the
raw RBF feature by
\[
\Phi_{ij}
=
\exp\left(
-\frac{\|z_i-c_j\|_2^2}{2h_j^2}
\right),
\qquad
i=1,\ldots,n,
\quad
j=1,\ldots,n_c.
\]

As in the one-dimensional setting, the feature columns are centered to
obtain \(F\), and the resulting test-function values are
\[
f = Fa.
\]

The derivative of the raw RBF feature with respect to the center is now
a two-dimensional vector:
\[
\nabla_{c_j}\Phi_{ij}
=
\Phi_{ij}\frac{z_i-c_j}{h_j^2}.
\]
The bandwidth derivative is
\[
\frac{\partial \Phi_{ij}}{\partial h_j}
=
\Phi_{ij}
\frac{\|z_i-c_j\|_2^2}{h_j^3}.
\]

Since centering is linear, the corresponding derivatives of the centered
features are obtained by applying the same centering operation to the
raw derivatives. Let
\[
r
=
\frac{\partial J}{\partial f}
=
\frac{1}{n}g
-
\frac{\|g\|_2^2}{n^2}f.
\]
Then the parameter gradients have the same chain-rule form as in the
one-dimensional case:
\[
\nabla_{c_j}J
=
a_j
\sum_{i=1}^{n}
r_i \nabla_{c_j}F_{ij},
\]
and
\[
\frac{\partial J}{\partial h_j}
=
a_j
\sum_{i=1}^{n}
r_i
\frac{\partial F_{ij}}{\partial h_j}.
\]

The centers and bandwidths are updated using the same gradient-ascent
and Armijo accept--reject procedure introduced in the one-dimensional
setting. The only difference is that each center update is now a vector
update in \(\mathbb{R}^2\):
\[
c_j
\leftarrow
c_j + \eta_c \nabla_{c_j}J.
\]

Thus, the two-dimensional formulation requires no change to the overall
optimization structure. The same construction extends directly to
\(z \in \mathbb{R}^d\) by replacing the two-dimensional vectors with
their \(d\)-dimensional counterparts.

\subsection{A two-dimensional nonlinear example}
\label{sec:two-dimensional-nonlinear-example}

In this part, we will test whether the complete nonlinear
basis construction remains effective for a genuinely two-dimensional
dependence pattern.  We generate
\begin{equation}
    z_i=(z_{i1},z_{i2})^{\mathsf T}
    \sim\operatorname{Unif}([-2,2]^2),
    \qquad
    x_i=s(z_i)+0.1\varepsilon_i,
    \qquad
    \varepsilon_i\sim\mathcal N(0,1),
    \label{eq:2d-two-bump-data}
\end{equation}
with the two-bump signal
\begin{equation}
    s(z_1,z_2)
    =\exp\!\left[-\frac{(z_1+0.7)^2+(z_2-0.4)^2}{0.25}\right]
     +0.7\exp\!\left[-\frac{(z_1-0.8)^2+(z_2+0.6)^2}{0.35}\right].
    \label{eq:2d-two-bump-signal}
\end{equation}
It therefore evaluates the ability of a multiple-center representation to adapt to a nontrivial two-dimensional
dependence function.

We use $n_c=9$ initial RBF centers arranged on the $3\times3$ Cartesian grid
\begin{equation}
    \mathcal C^{(0)}
    =\{-1.2,0,1.2\}\times\{-1.2,0,1.2\}
    \subset[-2,2]^2.
    \label{eq:2d-initial-center-grid}
\end{equation}
Thus, the nearest-neighbor spacing of the initial dictionary is
$d_{\rm NN}=1.2$.  All nine RBFs are initialized with the same isotropic
bandwidth,
\begin{equation}
    h_j^{(0)}=h_0=0.75,
    \qquad j=1,\ldots,9.
    \label{eq:2d-initial-bandwidth}
\end{equation}
This relatively broad initialization provides substantial overlap
between adjacent features and avoids gaps in the coverage of the
two-dimensional domain.

Starting from this common space-filling dictionary, we compare a fixed RBF
basis with the alternating center--bandwidth update.  In the fixed method,
the nine centers and all bandwidths remain at their initial values.  In the
adaptive method, each center is updated as a vector in $\mathbb R^2$ and each
RBF acquires its own scalar isotropic bandwidth $h_j$.  The bandwidth update
is performed through $q_j=\log h_j$ and is constrained to
\begin{equation}
    0.12\leq h_j\leq1.50,
    \label{eq:2d-bandwidth-bounds}
\end{equation}
while the centers are kept inside $[-2,2]^2$.  These bounds are numerical
safeguards rather than additional model assumptions.

\begin{table}[H]
    \centering
    \caption{Fixed and updated RBF bases for the two-dimensional two-bump
    example. Both methods use the same nine initial features and three
    detected dependence directions.}
    \label{tab:2d-two-bump-results}
    \begin{tabular}{lrrrr}
        \toprule
        Method
        & $J_{\rm Step\,1}$
        & $|\operatorname{corr}(y,s)|$
        & $\operatorname{corr}(y,\varepsilon)$
        & Transport cost \\
        \midrule
        Fixed basis
        & 0.1770 & 0.4910 & 0.7229 & 0.00919 \\
        Updated basis
        & 0.4175 & 0.0262 & 0.9882 & 0.01683 \\
        \bottomrule
    \end{tabular}
\end{table}

\begin{figure}[H]
    \centering
    \includegraphics[width=\textwidth]
        {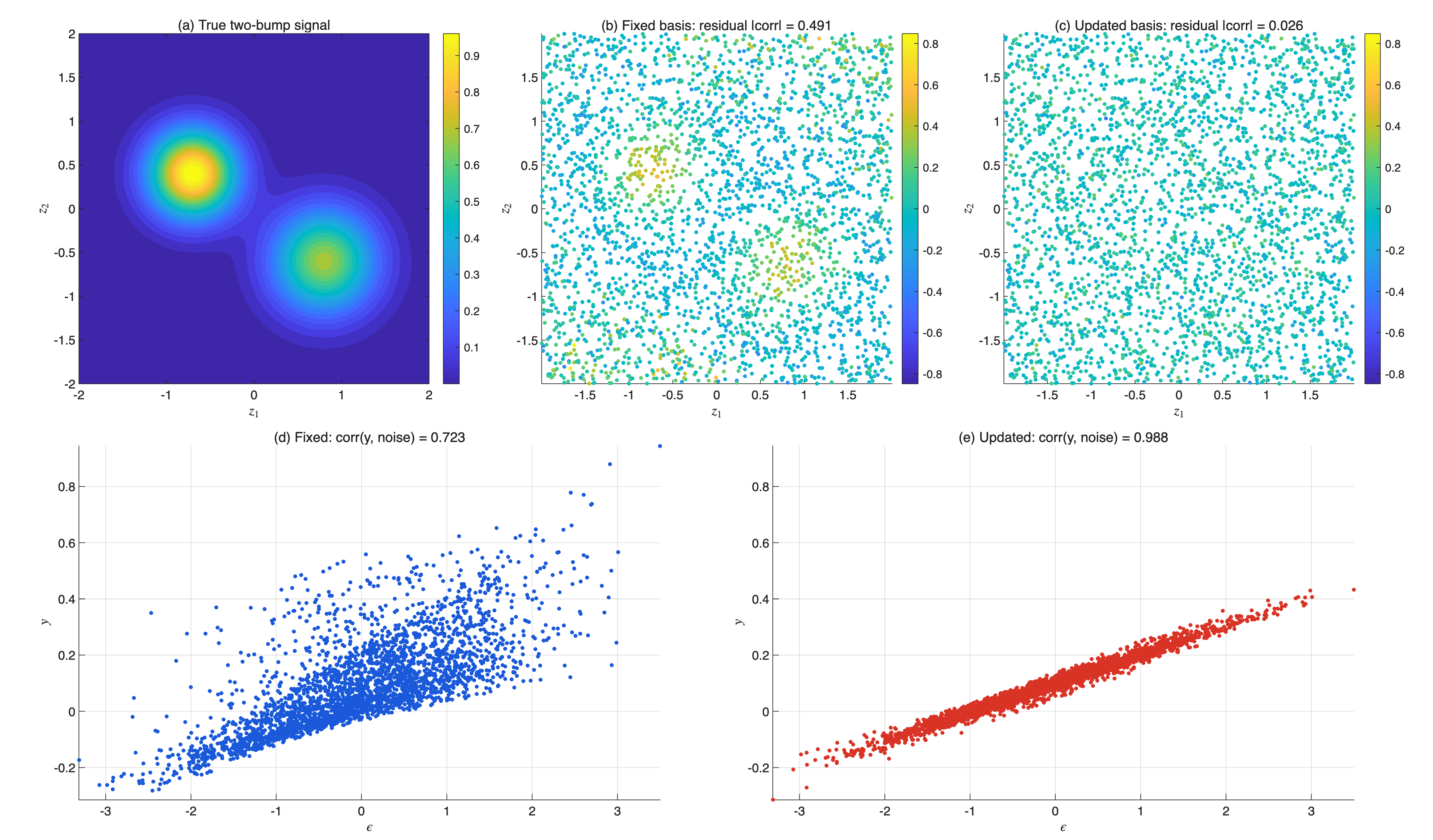}
    \caption{Two-dimensional two-bump example.  The top row shows the true
    signal and the centered transported values over the $(z_1,z_2)$-plane
    for the fixed and adaptive RBF bases, using a common color scale for the
    two transported-value panels.  The bottom row compares the transported
    variable with the independent noise.  The adaptive basis removes the
    localized two-dimensional pattern and produces a transported variable
    that is nearly a function of the independent noise alone.  The associated
    center and bandwidth changes are shown separately in
    Figure~\ref{fig:2d-basis-update}.}
    \label{fig:2d-two-bump-comparison}
\end{figure}

\begin{figure}[H]
    \centering
    \includegraphics[width=\textwidth]
        {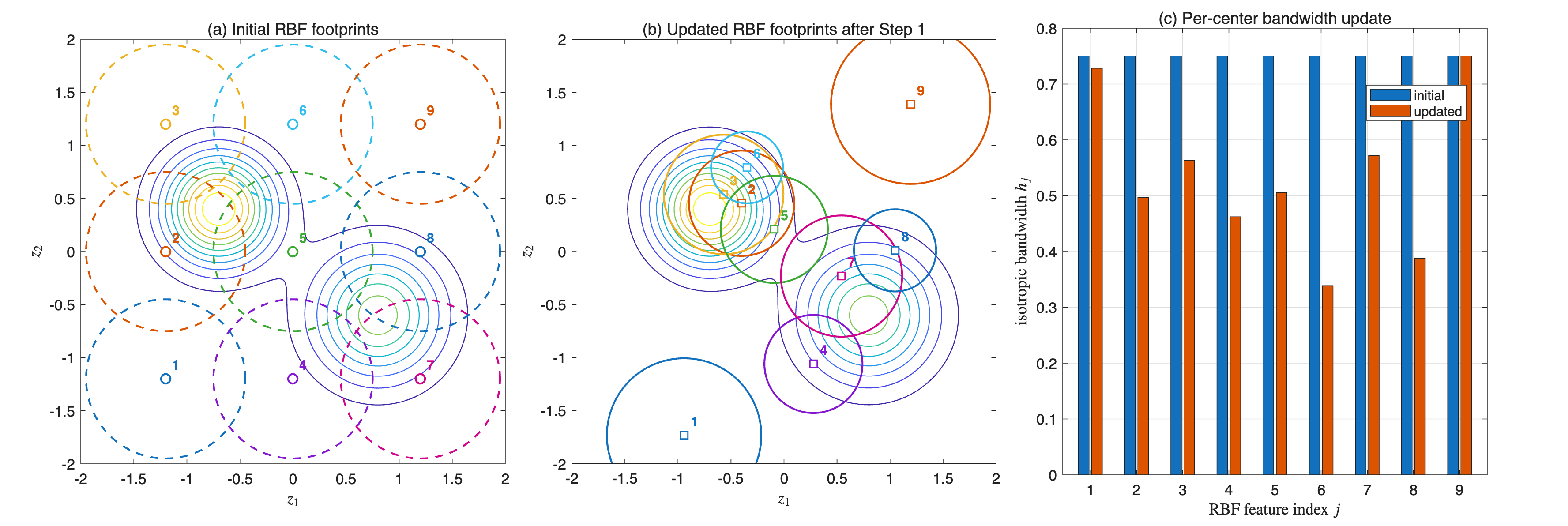}
    \caption{Change in the adaptive RBF dictionary during the first Step~1
    dependence search.  Left: initial centers and isotropic RBF footprints.
    Middle: updated centers and footprints, with each circle showing the
    $\exp(-1/2)$ level set and hence having radius $h_j$.  Center labels
    identify the same feature in both panels.  Right: initial and updated
    per-center bandwidths.  The separate spatial panels show where the
    features move, while the bar plot makes the bandwidth changes directly
    comparable.}
    \label{fig:2d-basis-update}
\end{figure}

As shown in Figure~\ref{fig:2d-basis-update}, centers 1 and 9
remain relatively far from the two bumps. Nevertheless, the
updated basis substantially reduces the residual dependence,
suggesting that accurate dependence removal does not require
every RBF center to move toward a signal peak.

Here, we deliberately begin with isotropic bandwidths.  This isolates the effect of
increasing the dimension of the conditioning variable while adding only one
new center coordinate.  Once this case is validated, an axis-aligned
anisotropic extension follows by replacing $h_j$ with
$(h_{j1},h_{j2})$ and using
\begin{equation}
    \Phi_{ij}
    =\exp\!\left[-\frac12\sum_{k=1}^2
        \frac{(z_{ik}-c_{jk})^2}{h_{jk}^2}\right].
    \label{eq:2d-diagonal-rbf}
\end{equation}
Its coordinate-wise derivatives are
\begin{equation}
    \frac{\partial\Phi_{ij}}{\partial c_{jk}}
    =\Phi_{ij}\frac{z_{ik}-c_{jk}}{h_{jk}^2},
    \qquad
    \frac{\partial\Phi_{ij}}{\partial h_{jk}}
    =\Phi_{ij}\frac{(z_{ik}-c_{jk})^2}{h_{jk}^3}.
    \label{eq:2d-diagonal-rbf-gradients}
\end{equation}
The same centering, chain rule, and Armijo procedure then apply without any
other change to the algorithm.
\FloatBarrier
\section{Real-data experiments}
\label{sec:real-data-experiments}
Having evaluated the method on synthetic examples, we next investigate its performance on real-world datasets in this section.
\subsection{Age-adjusted alkaline phosphatase in NHANES}
\label{sec:nhanes-alp-age}

Age adjustment and standardization are routinely used in epidemiology to
compare populations whose age compositions differ
\citep{cdcAgeAdjustment2026}. In this example, We consider serum alkaline phosphatase (ALP) measurements from the National Health and Nutrition Examination Survey (NHANES) as our first real-data example. ALP is an enzyme primarily associated with the liver and bones, and is commonly measured in routine blood tests. ALP levels are typically higher in children and adolescents due to rapid bone growth and development. As skeletal growth slows and stabilizes in adulthood, ALP levels generally decrease. As a result, the distribution of ALP measurements can differ substantially across age groups. Our goal is to construct an age-adjusted variable (y=T(x,z)) that is approximately independent of age (z), while remaining as close as possible to the original ALP measurements (x). This allows us to examine whether OTBP can reduce age-related differences in the conditional distribution of a real biomarker.

\subsubsection{Data and preprocessing}
\label{sec:nhanes-alp-data}

We use the public 2017--2018 National Health and Nutrition Examination Survey
(NHANES).  Serum alkaline phosphatase (ALP) is taken from the Standard
Biochemistry Profile, file \texttt{BIOPRO\_J}, where \texttt{LBXSAPSI}
records ALP in IU/L \citep{nhanes2017biopro}.  Age and demographic variables
are obtained from \texttt{DEMO\_J}, and the two files are merged by the
participant identifier \texttt{SEQN}.  The eligible laboratory sample begins
at age 12.  We retain participants aged 12--79 with finite positive ALP and
exclude participants recorded as pregnant.  Age 80 and above is omitted
because NHANES top-codes those ages.  From the remaining observations, we
draw a reproducible random sample of $n=3000$ participants.

Because ALP is positive and right-skewed, the transported outcome is the
standardized log measurement
\begin{equation}
    x_i=\frac{\log(\mathrm{ALP}_i)-\overline{\ell}}{s_{\ell}},
    \qquad \ell_i=\log(\mathrm{ALP}_i),
    \label{eq:nhanes-alp-outcome}
\end{equation}
where $\overline{\ell}$ and $s_{\ell}$ are the sample mean and standard
deviation of log ALP.  The conditioning varbable is standardized age,
\begin{equation}
    z_i=\frac{\mathrm{Age}_i-\overline{\mathrm{Age}}}
    {s_{\mathrm{Age}}}.
    \label{eq:nhanes-age-covariate}
\end{equation}
The analysis is intentionally unweighted and is used as a methodological
sample-level illustration.  It is not a population estimate, a clinical
reference-range analysis, or a diagnostic procedure.

\subsubsection{OTBP specification}
\label{sec:nhanes-alp-specification}

Since age is one-dimensional, we use a fixed dictionary of $n_c=9$ Gaussian
RBFs.  Their centers are equally spaced over the standardized-age range,
\begin{equation}
    c_j\in\{-2.4,-1.8,-1.2,-0.6,0,0.6,1.2,1.8,2.4\},
    \qquad h_j=0.60,
    \label{eq:nhanes-alp-rbf}
\end{equation}
for $j=1,\ldots,9$.  Preliminary center and bandwidth optimization produced
nearly identical dependence diagnostics, indicating that this dense fixed
dictionary was already sufficiently expressive; we therefore report the fixed-basis result only for $n_c=9$.  The $G$-feature family is
\begin{equation}
    G(y)=\bigl(y,y^2,y^3,y^4\bigr),
\end{equation}
and the sequential search uses $K_{\max}=8$ as an upper bound, with early
stopping under the standard detection rule.  Three dependence directions are
retained in this experiment.

\subsubsection{Evaluation criteria} 
\label{sec:real-data-evaluation}  
For both real-data examples, we report three complementary diagnostics of residual dependence between the transported variable and the conditioning covariate. First, distance correlation measures general nonlinear association without requiring a parametric model.\footnote{
Given observations $\{(y_i,z_i)\}_{i=1}^n$, define
$a_{ij}=\|z_i-z_j\|_2$ and
$b_{ij}=|y_i-y_j|$.
Let $A$ and $B$ denote the corresponding
doubly centered distance matrices.
The sample distance covariance is defined by
\[
\operatorname{dCov}_n^2(Y,Z)
=\frac{1}{n^2}\sum_{i,j=1}^n A_{ij}B_{ij},
\]
and the sample distance correlation is
\[
\operatorname{dCor}_n(Y,Z)
=
\sqrt{
\frac{\operatorname{dCov}_n^2(Y,Z)}
{\operatorname{dCov}_n(Y,Y)
 \operatorname{dCov}_n(Z,Z)}
}.
\]
The population distance correlation is zero
if and only if $Y$ and $Z$ are independent,
provided the relevant first moments are finite.
} Second, normalized HSIC uses Gaussian kernels with bandwidths selected by the median-distance rule; smaller values indicate weaker residual dependence. Third, a permutation diagnostic compares the observed normalized HSIC with 100 random permutations. All diagnostics are computed on the same random subset of $m=800$ observations for each method. With 100 permutations, the smallest attainable reported $p$-value is $1/101\simeq0.0099$.  We also report the transport cost \begin{equation}     \mathcal C(x,y)=\frac{1}{2n}\lVert y-x\rVert_2^2,     \label{eq:real-data-transport-cost} \end{equation} the objective from the first Step~1 search, and the number of retained dependence directions. Conditional location and scale are assessed by partitioning the conditioning covariate into bins and comparing the within-bin means and standard deviations of the original and transported outcomes. For the ALP analysis, age is divided into ten sample-quantile bins.

\subsubsection{Results}
\label{sec:nhanes-alp-results}

Table~\ref{tab:nhanes-alp-results} summarizes the dependence diagnostics.

\begin{table}[H]
    \centering
    \caption{NHANES ALP--age results.  Dependence diagnostics are computed on
    the common $m=800$ subset, using 100 permutations for HSIC.}
    \label{tab:nhanes-alp-results}
        \resizebox{\textwidth}{!}{%

    \begin{tabular}{lrrrrrr}
        \toprule
        Method & First Step~1 $J$ & Distance corr. & Normalized HSIC
        & Permutation $p$ & Transport cost & $K_{\rm used}$ \\
        \midrule
        Original $x$ & --     & 0.2879 & 0.05157  & 0.0099 & 0      & -- \\
        Fixed basis, $n_c=9$ & 0.2661 & 0.0375 & 0.000512 & 0.9901 & 0.1925 & 3 \\
        Fixed basis, $n_c=3$ & 0.1176 & 0.1494 & 0.01914 & 0.0099 & 0.1053 & 3 \\
        Updated basis, $n_c=3$ & 0.2686 & 0.0387 & 0.000777 & 0.9604 & 0.1996 & 4 \\
        \bottomrule
    \end{tabular}%
    }\end{table}

With $n_c=9$ RBF centers, the fixed basis already achieves substantial
dependence reduction, yielding a distance correlation of $0.0375$ and a
normalized HSIC of $0.000512$. To investigate whether basis optimization
can achieve comparable performance with fewer features, we reduce the
number of centers to $n_c=3$. In this smaller feature space, the fixed
basis performs considerably worse, with distance correlation increasing
to $0.1494$. However, updating the RBF centers and bandwidths reduces
the distance correlation to $0.0387$ and the normalized HSIC to $0.000777$,
achieving dependence removal comparable to that of the nine-center fixed
basis. These results suggest that basis optimization can compensate for
a reduced number of features by improving their ability to capture the
underlying dependence structure.

\begin{figure}[H]
    \centering
    \includegraphics[width=0.94\textwidth]
        {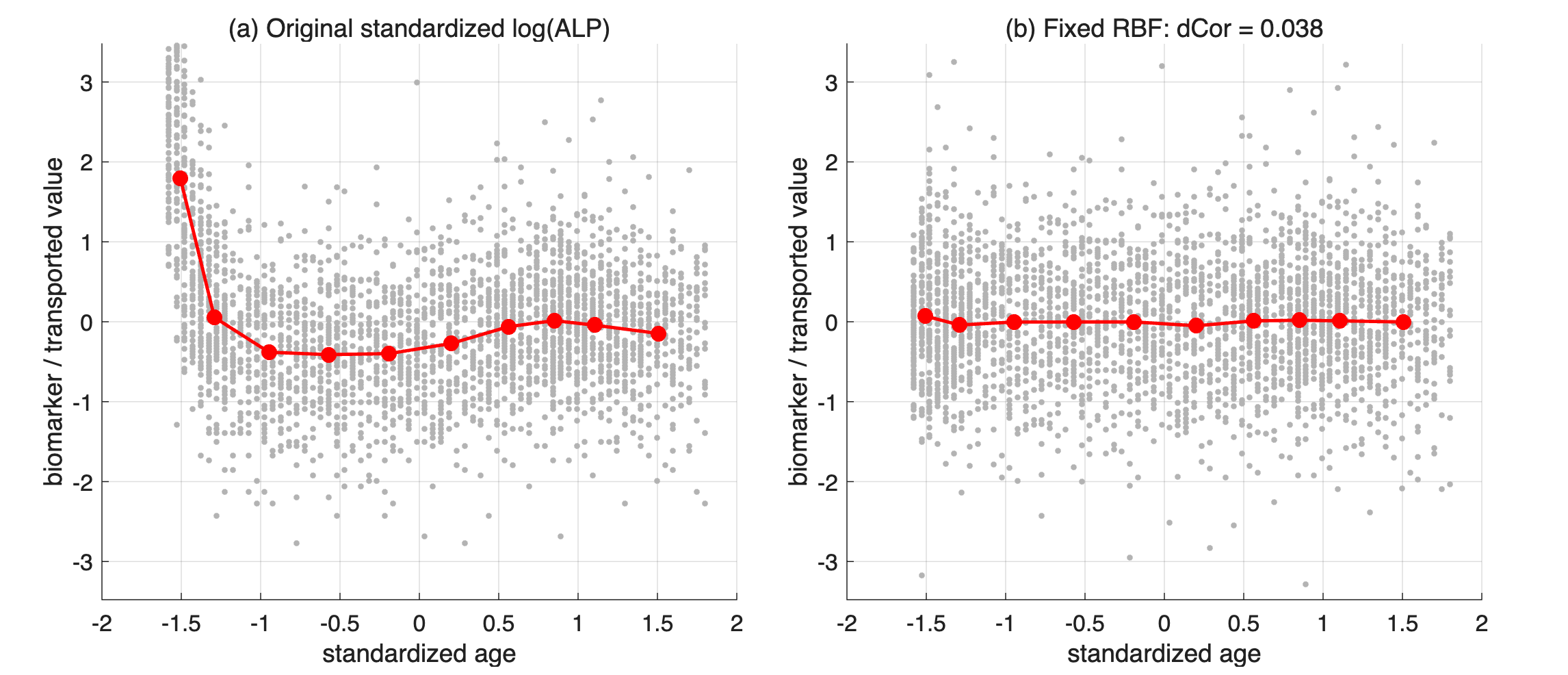}
    \caption{Original and transported standardized log ALP over standardized
    age.  Red points and lines show age-bin means.  Left: original outcome.
    Right: fixed-RBF OTBP transport.}
    \label{fig:nhanes-alp-dependence}
\end{figure}

\begin{figure}[H]
    \centering
    \includegraphics[width=0.82\textwidth]
        {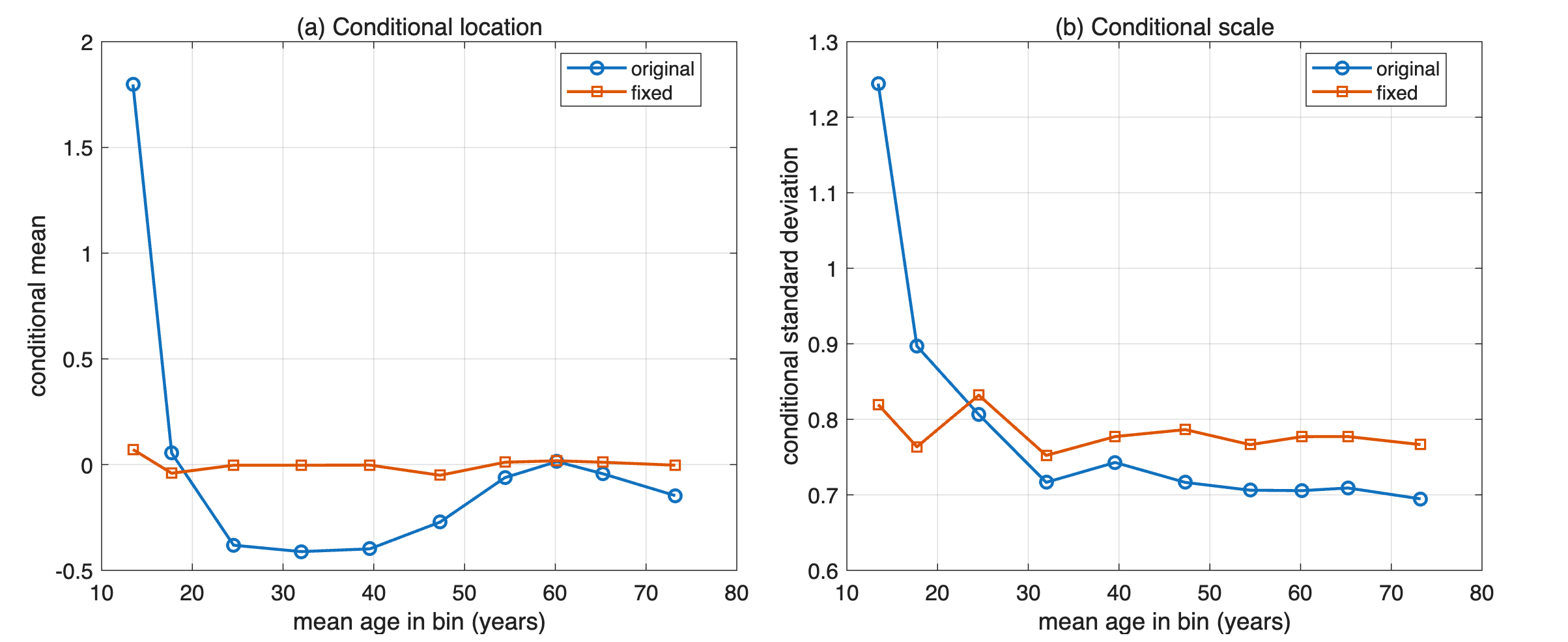}
    \caption{Conditional location and scale in the NHANES ALP example.
    Points are computed in ten sample-quantile age bins.  The transport
    strongly attenuates the age-dependent conditional mean and compresses
    the between-bin variation in conditional standard deviation.}
    \label{fig:nhanes-alp-conditional}
\end{figure}\begin{figure}[H]
    \centering
    \includegraphics[width=0.98\textwidth]
        {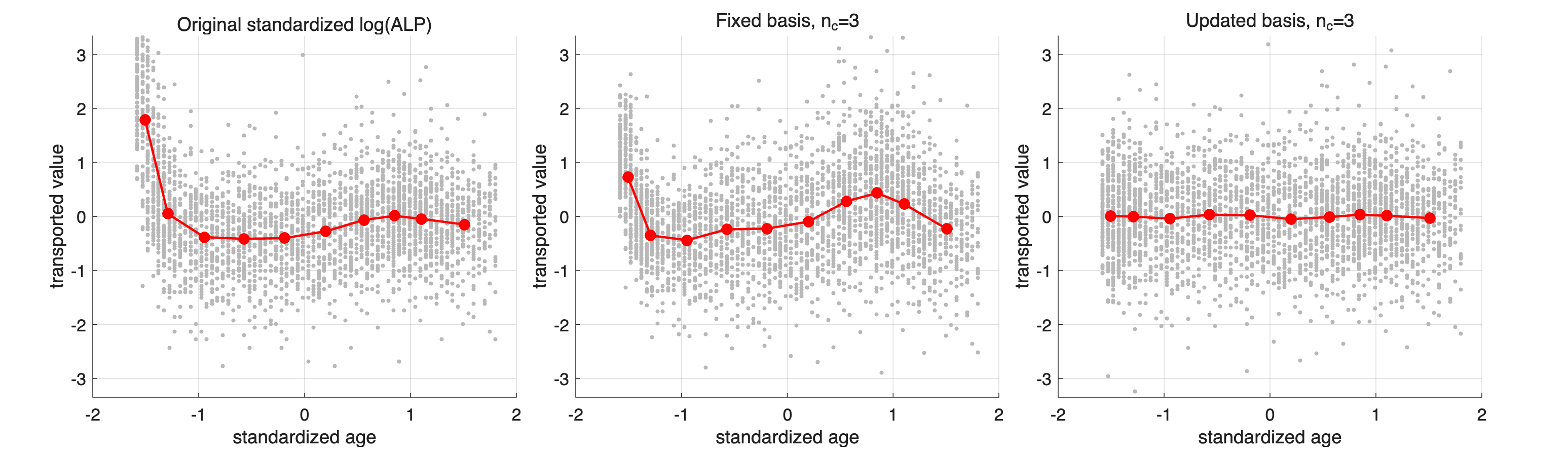}
    \caption{Effect of restricting the NHANES ALP dictionary to three RBF centers. Red points and lines show means in ten age-quantile bins. Left: original standardized log ALP. Middle: transport with the fixed three-center basis. Right: transport after updating the three centers and their bandwidths.}
    \label{fig:nhanes-alp-nc3}
\end{figure}

\begin{figure}[H]
    \centering
    \includegraphics[width=0.98\textwidth]
        {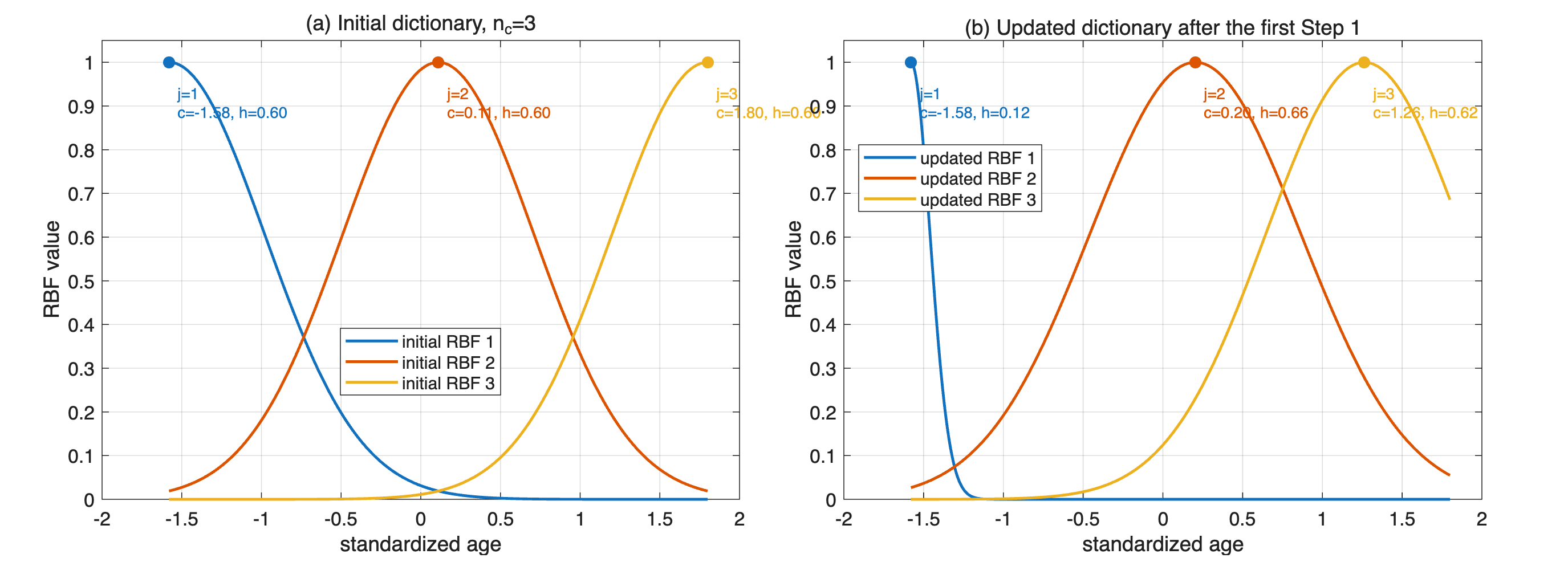}
    \caption{Update of the three-center RBF dictionary in the NHANES ALP example. Left: the initial Gaussian RBF features. Right: the centers and bandwidths after the first Step~1 dependence search. The annotations report the center and bandwidth of each feature in standardized age coordinates.}
    \label{fig:nhanes-alp-nc3-basis-update}
\end{figure}

The corresponding center and bandwidth changes for the adaptive three-center dictionary are shown in Figure~\ref{fig:nhanes-alp-nc3-basis-update}.

\subsection{Bike-sharing demand}
\label{sec:bike-sharing}

We use the hourly Bike Sharing data from the UCI Machine Learning Repository \citep{fanaee2013bike}.  The data contain the total number of casual and registered bicycle rentals together with weather and calendar variables for 2011--2012. We test whether our methods can remove the empirical dependence of hourly rental demand on a genuinely two-dimensional environmental covariate.

\subsubsection{Data and preprocessing}
\label{sec:bike-sharing-data}

Let $\mathrm{cnt}_i$ denote the total number of rentals in hour $i$.  Because
the raw counts are strongly right-skewed, we use the standardized log count
as the transported variable,
\begin{equation}
    x_i
    =\frac{\log(1+\mathrm{cnt}_i)-\overline{\ell}}{s_{\ell}},
    \qquad
    \ell_i=\log(1+\mathrm{cnt}_i),
    \label{eq:bike-outcome}
\end{equation}
where $\overline{\ell}$ and $s_{\ell}$ are the sample mean and standard
deviation of the log counts.  The conditioning variable is
\begin{equation}
    z_i=(z_{i1},z_{i2})^{\mathsf T},
    \label{eq:bike-covariate}
\end{equation}
where $z_{i1}$ and $z_{i2}$ are the standardized temperature and humidity,
respectively.  Although these variables are supplied in normalized form in
the original data, we standardize each coordinate again within the selected
sample so that a common isotropic RBF bandwidth has the same interpretation
in both directions.

To retain the within-day structure of the hourly observations, complete
calendar dates rather than isolated hours are sampled as blocks.  Dates are
randomly added until the target sample size of 3000 is reached, resulting in
$n=3023$ hourly observations.  The same sampled dates, preprocessing
constants, initial dictionary, and random seed are used for every method in
the comparison.

\subsubsection{OTBP specification}
\label{sec:bike-sharing-specification}

We compare the original standardized log count with fixed RBF basis and updated RBF basis method.

Both variants start from $n_c=9$ isotropic RBFs whose centers form the
$3\times3$ grid
\begin{equation}
    \mathcal C^{(0)}
    =\{-1.2,0,1.2\}\times\{-1.2,0,1.2\},
    \label{eq:bike-initial-centers}
\end{equation}
in the standardized temperature--humidity plane.  Every feature is
initialized at $h_j^{(0)}=0.75$.  During adaptive optimization the centers
are constrained to $[-3,3]^2$ and the isotropic bandwidths are constrained
to
\begin{equation}
    0.12\leq h_j\leq1.50.
    \label{eq:bike-bandwidth-bounds}
\end{equation}
For each detected direction, the adaptive method uses four alternating
blocks, with at most 150 accepted center proposals followed by 150 bandwidth
proposals in each block.  The same Armijo accept--reject rule used in the
simulation experiments is applied here.

To detect dependence involving more than the first two conditional moments,
the $G$-feature family contains polynomial terms through degree four,
\begin{equation}
    G(y)=\bigl(y,y^2,y^3,y^4\bigr).
    \label{eq:bike-G-features}
\end{equation}
We set $K_{\max}=8$, where $K_{\max}$ is only an upper bound on the number of
sequentially detected dependence directions.  The search terminates early
when the next direction falls below the same detection threshold used by the
OTBP algorithm.  Consequently, increasing $K_{\max}$ does not force the
method to use eight directions.

\subsubsection{Evaluation criteria}
\label{sec:bike-sharing-evaluation}

We use the same distance-correlation, normalized-HSIC, permutation, and
transport-cost criteria defined in Section~\ref{sec:real-data-evaluation},
with the same common diagnostic subset of $m=800$ observations and 100 HSIC
permutations. The only example-specific summary is the conditional
location--scale display: temperature and humidity are each divided into
sample quartiles, producing a $4\times4$ grid of 16 cells. Within each cell,
we compute the mean and standard deviation of the original and transported
outcomes. Approximate independence should make these summaries relatively
stable across cells, up to sampling variability.

\subsubsection{Results}
\label{sec:bike-sharing-results}

\begin{table}[H]
    \centering
    \caption{Bike Sharing results for the original outcome and the fixed and
    adaptive RBF transports.  The Step~1 objective is recorded for the first
    detected direction.  The permutation diagnostic uses 100 permutations
    of the common $m=800$ diagnostic subset.}
    \label{tab:bike-main-results}
    \resizebox{\textwidth}{!}{%
    \begin{tabular}{lrrrrrr}
        \toprule
        Method
        & First Step~1 $J$
        & Distance corr.
        & Normalized HSIC
        & Permutation $p$
        & Transport cost
        & $K_{\rm used}$ \\
        \midrule
        Original $x$
        & --     & 0.4987 & 0.22046 & 0.0099 & 0      & -- \\
        Fixed RBF basis
        & 0.1588 & 0.1459 & 0.01333 & 0.0099 & 0.1510 & 5 \\
        Updated RBF basis
        & 0.1813 & 0.0929 & 0.00613 & 0.1386 & 0.1765 & 7 \\
        \bottomrule
    \end{tabular}%
    }
\end{table}

The fixed-basis residual remains significant under the permutation
diagnostic, whereas the updated-basis value is not significant at the
$5\%$ level.  

\begin{figure}[H]
    \centering
    \includegraphics[width=\textwidth]
        {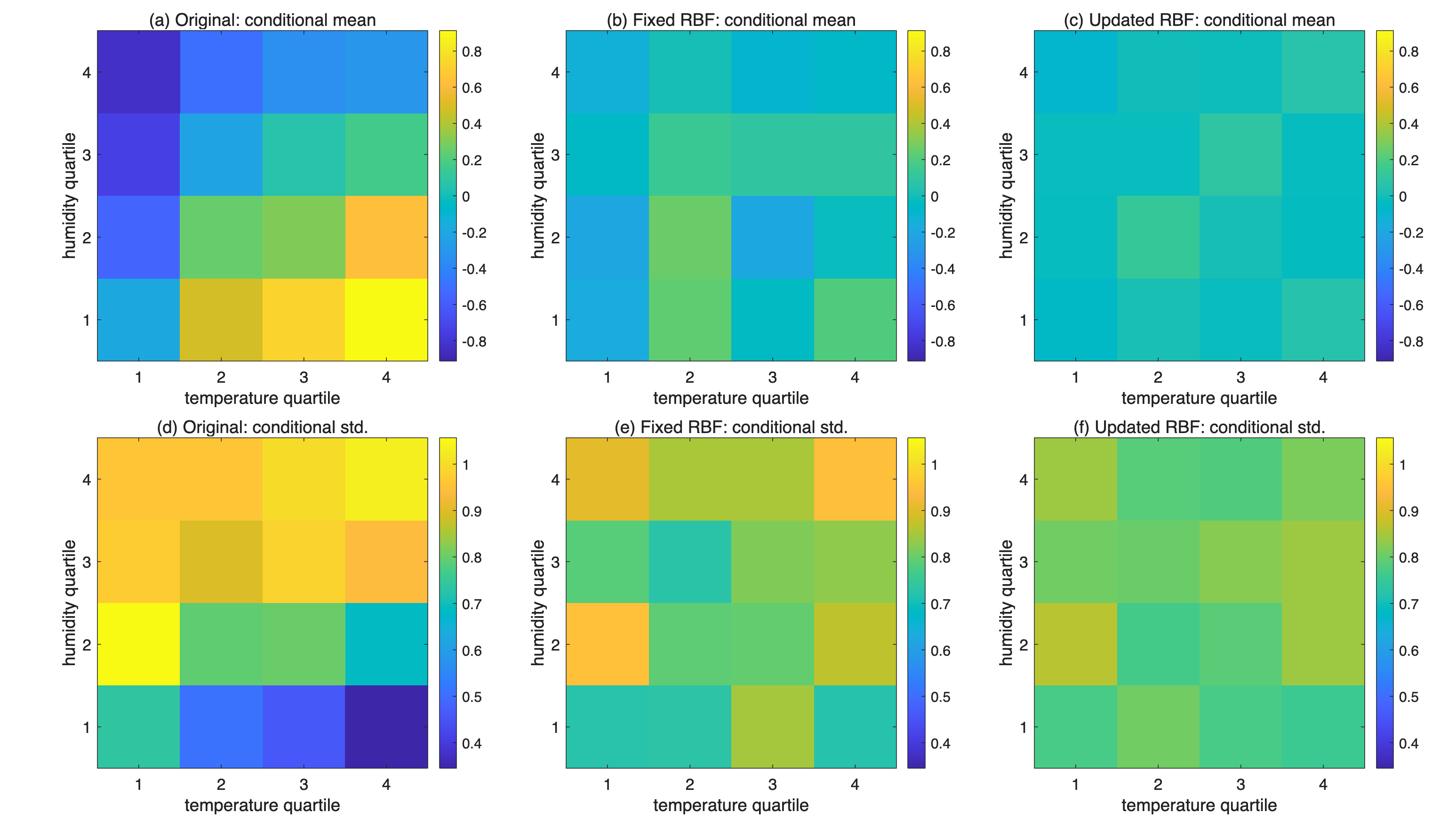}
    \caption{Conditional location and scale in the Bike Sharing example.
    Temperature and humidity are each divided into sample quartiles.  The top
    row shows cell-wise conditional means, and the bottom row shows cell-wise
    conditional standard deviations, for the original outcome, fixed-basis
    transport, and adaptive-basis transport.  A common color scale is used
    within each row.}
    \label{fig:bike-conditional-statistics}
\end{figure}

\begin{figure}[H]
    \centering
    \includegraphics[width=0.98\textwidth]
        {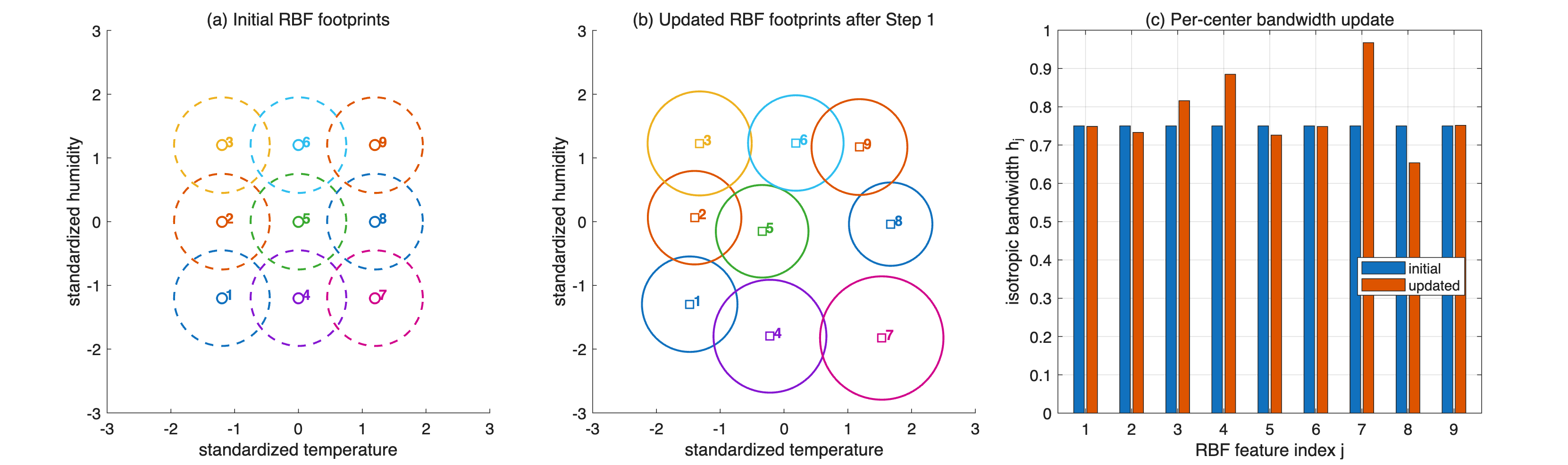}
    \caption{Adaptive RBF dictionary update in the Bike Sharing example. Left: initial centers and isotropic RBF footprints in standardized temperature--humidity coordinates. Middle: updated centers and footprints after the first Step~1 search. Right: initial and updated per-center isotropic bandwidths. Matching labels identify the same feature across panels.}
    \label{fig:bike-sharing-basis-update}
\end{figure}

Figure~\ref{fig:bike-sharing-basis-update} shows how the centers and isotropic bandwidths change during the first Step~1 dependence search.

Overall, the real-data example supports the same conclusion as the
two-dimensional simulation: adapting the centers and local scales improves
the ability of a finite RBF dictionary to identify and remove nonlinear
dependence. 

\section{Discussion}

A central point of this work is that the independence constraint in the
sample-based OTBP can be estimated through a finite-dimensional approximation of the test-function spaces. Rather than testing independence over all admissible functions, the numerical method searches for dependence within chosen feature spaces,
\[
f=Fa,
\qquad
g=Gb.
\]
The choice of these feature spaces therefore determines which forms of
dependence can be detected and removed. This feature-space approximation is
the key link between the theoretical independence condition and the practical finite-sample algorithm.

The linear setting provides a simple and interpretable realization of this
idea, but its ability to detect dependence is limited by the selected
low-order features. Replacing the linear feature matrix for $f$ by Gaussian
RBF features substantially enlarges the detectable class without changing the underlying OTBP formulation. The numerical examples show that the nonlinear feature representation can capture oscillatory, localized, monotone, and nonlinearly heteroscedastic dependence and so on.

Another contribution is to allow the nonlinear basis itself to adapt to the Step~1 dependence objective. Rather than keeping the RBF centers and
bandwidths fixed after initialization, we update the centers and per-center
bandwidths by gradient ascent on the dependence objective, using an
Armijo-type accept--reject rule \citep{armijo1966minimization} and a block-wise alternating scheme. They show that basis optimization can compensate for a relatively small feature space: in the sinusoidal example, an updated basis with fewer centers reaches a low-residual regime comparable to that of a substantially larger fixed basis.

The experiments with non-uniform covariates also highlight the distinction
between density-aware and dependence-aware feature placement. In the
localized middle-bump example with a mixture-Gaussian distribution of $z$,
uniform initialization performs better than $k$-means after the basis updates. This suggests that a density-based initialization may place too much emphasis on high-density regions, even though lower-density regions may still contain important dependence. A broader initial coverage can give the subsequent dependence-driven updates more flexibility to move the basis toward locations that are useful for detecting and removing dependence.

We further extended the adaptive RBF formulation to multidimensional covariates, demonstrating its effectiveness in a two-dimensional synthetic example with the same optimization framework applicable to higher dimensions. Real-data experiments involving age-adjusted biomarker measurements and bike-sharing demand further illustrate the applicability of the method. These results suggest that basis adaptation can improve dependence removal when the initial feature space is insufficiently expressive, though sometimes at the expense of increased transportation cost.

Several limitations remain. While this work develops adaptive nonlinear features for $f(z)$, the representation of $g(y)$ is restricted to fixed low-order polynomial features, which may limit the detection of more complex dependence structures. Moreover, the basis optimization is
nonconvex and sensitive to initialization, and formal convergence guarantees for the alternating center--bandwidth updates remain to be established. Although the formulation extends to multidimensional covariates, the numerical experiments are limited to only isotropic RBFs, with a single bandwidth shared across all coordinate directions for each basis function. 
\par Future work should investigate more flexible feature representations for higher-dimensional covariates, richer $G$-feature families, theoretical convergence properties, repeated real-data analyses with uncertainty estimates, and formal out-of-sample evaluation.

\bibliographystyle{elsarticle-harv}
\bibliography{references}

@article{agueh2011barycenter,
  author  = {Agueh, Martial and Carlier, Guillaume},
  title   = {Barycenters in the {W}asserstein Space},
  journal = {SIAM Journal on Mathematical Analysis},
  volume  = {43},
  number  = {2},
  pages   = {904--924},
  year    = {2011},
  doi     = {10.1137/100805741}
}

@article{pass2013infinitely,
  author  = {Pass, Brendan},
  title   = {Optimal Transportation with Infinitely Many Marginals},
  journal = {Journal of Functional Analysis},
  volume  = {264},
  number  = {4},
  pages   = {947--963},
  year    = {2013},
  doi     = {10.1016/j.jfa.2012.12.002}
}

@article{tabak2018variability,
  author  = {Tabak, Esteban G. and Trigila, Giulio},
  title   = {Explanation of Variability and Removal of Confounding Factors from Data through Optimal Transport},
  journal = {Communications on Pure and Applied Mathematics},
  volume  = {71},
  number  = {1},
  pages   = {163--199},
  year    = {2018},
  doi     = {10.1002/cpa.21706}
}

@article{tabak2022distributional,
  author  = {Tabak, Esteban G. and Trigila, Giulio and Zhao, Wenjun},
  title   = {Distributional Barycenter Problem through Data-Driven Flows},
  journal = {Pattern Recognition},
  volume  = {130},
  pages   = {108795},
  year    = {2022},
  doi     = {10.1016/j.patcog.2022.108795}
}

@article{lipnick2026monge,
  author  = {Lipnick, Andrew D. and Tabak, Esteban G. and Trigila, Giulio and Wang, Yating and Ye, Xuancheng and Zhao, Wenjun},
  title   = {The {M}onge Optimal Transport Barycenter Problem},
  journal = {International Statistical Review},
  year    = {2026},
  note    = {To appear}
}

@book{silverman1986density,
  author    = {Silverman, Bernard W.},
  title     = {Density Estimation for Statistics and Data Analysis},
  publisher = {Chapman and Hall},
  year      = {1986}
}

@article{armijo1966minimization,
  author  = {Armijo, Larry},
  title   = {Minimization of Functions Having {Lipschitz} Continuous First Partial Derivatives},
  journal = {Pacific Journal of Mathematics},
  volume  = {16},
  number  = {1},
  pages   = {1--3},
  year    = {1966},
  doi     = {10.2140/pjm.1966.16.1}
}

@article{lloyd1982least,
  author  = {Lloyd, Stuart P.},
  title   = {Least Squares Quantization in {PCM}},
  journal = {IEEE Transactions on Information Theory},
  volume  = {28},
  number  = {2},
  pages   = {129--137},
  year    = {1982},
  doi     = {10.1109/TIT.1982.1056489}
}

@article{micchelli1986interpolation,
  author  = {Micchelli, Charles A.},
  title   = {Interpolation of Scattered Data: Distance Matrices and Conditionally Positive Definite Functions},
  journal = {Constructive Approximation},
  volume  = {2},
  number  = {1},
  pages   = {11--22},
  year    = {1986},
  doi     = {10.1007/BF01893414}
}

@inproceedings{gretton2005measuring,
  author    = {Gretton, Arthur and Bousquet, Olivier and Smola, Alexander and Sch{\"o}lkopf, Bernhard},
  title     = {Measuring Statistical Dependence with {Hilbert--Schmidt} Norms},
  booktitle = {Algorithmic Learning Theory},
  series    = {Lecture Notes in Computer Science},
  volume    = {3734},
  pages     = {63--77},
  publisher = {Springer},
  year      = {2005},
  doi       = {10.1007/11564089_7}
}

@misc{fanaee2013bike,
  author    = {Fanaee-T, Hadi},
  title     = {Bike Sharing},
  year      = {2013},
  howpublished = {UCI Machine Learning Repository [Dataset]},
  doi       = {10.24432/C5W894},
  url       = {https://doi.org/10.24432/C5W894}
}

@misc{cdcAgeAdjustment2026,
  author       = {{National Center for Health Statistics}},
  title        = {Age Adjustment},
  year         = {2026},
  organization = {Centers for Disease Control and Prevention},
  url          = {https://www.cdc.gov/nchs/hus/sources-definitions/age-adjustment.htm},
  note         = {Accessed 2026-10-29}
}

@misc{nhanes2017biopro,
  author       = {{National Center for Health Statistics}},
  title        = {{NHANES} 2017--2018: Standard Biochemistry Profile Data Documentation, Codebook, and Frequencies},
  year         = {2020},
  organization = {Centers for Disease Control and Prevention},
  url          = {https://wwwn.cdc.gov/Nchs/Data/Nhanes/Public/2017/DataFiles/BIOPRO_J.htm},
  note         = {Accessed 2026-10-29}
}

\end{document}